\documentclass[final]{elsarticle} % switch final for review

\usepackage{amssymb}
\usepackage{amsmath, amsthm}
\usepackage{graphicx}
\usepackage{enumitem}
\usepackage{algorithm}
\usepackage[noend]{algpseudocode}
\usepackage{xcolor}
\usepackage{multirow}
\usepackage{eqparbox}
\usepackage{subcaption}
\usepackage{bm}
\usepackage{booktabs}
\usepackage{hyperref}

\newtheorem{theorem}{Theorem}
\newtheorem{definition}{Definition}

\usepackage[textsize=tiny]{todonotes}

\usepackage[letterpaper,centering]{geometry}
\newcommand{\revise}[1]{{#1}}

\newcolumntype{C}{>{\centering\arraybackslash} m{4cm} }  %# New column type
\newcolumntype{D}{>{\centering\arraybackslash} m{5cm} }  %# New column type  

\title{A continuous analogue of the tensor-train decomposition}
\author[myaddress]{Alex Gorodetsky\corref{mycor}}
\cortext[mycor]{Corresponding author}
\ead{goroda@umich.edu}

\author[mit]{Sertac Karaman}
\ead{sertac@mit.edu}
\author[mit]{Youssef Marzouk}
\ead{ymarz@mit.edu}

\address[myaddress]{University of Michigan, Ann Arbor, MI 48109, USA}
\address[mit]{Massachusetts Institute of Technology, Cambridge, MA 02139, USA}

\newcommand{\posint}{\mathbb{Z}^{+}}
\newcommand{\reals}{\mathbb{R}}
\newcommand{\xspace}[1]{\mathcal{X}_{#1}}

\newcommand{\mat}[1]{\mathbf{#1}}

\newcommand{\bvec}[1]{\bm{#1}}
\newcommand{\mvf}[1]{\mathcal{#1}} % matrix valued function

\newcommand{\submatm}[1]{\bar{\mat{#1}}_{\textnormal{mat}}} %matrix valued function

\newcommand{\foldleft}[1]{x_{\leq#1}}

\newcommand{\foldright}[1]{x_{>#1}}
\newcommand{\frspace}[1]{\xspace{>#1}}

\newcommand{\funfold}[1]{f^{#1}} % unfolding function

\newcommand{\fcore}[1]{\mvf{F}_{#1}(x_{#1})} % function train core
\newcommand{\facore}[1]{\hat{\mvf{F}}_{#1}(x_{#1})} % function train core
\newcommand{\ffiber}[2]{f_{#1}^{(#2)}} % fiber
\newcommand{\fafiber}[2]{\hat{f}_{#1}^{(#2)}} % fiber
\newcommand{\efiber}[2]{\Delta_{#1}^{(#2)}} % fiber
\newcommand{\fn}[1]{n_{#1}^{(\alpha_{#1-1}, \alpha_{#1})}}
\newcommand{\fnn}[1]{n_{#1}}
\newcommand{\parambase}{\theta}
\newcommand{\paramk}[2]{\parambase_{#1#2}^{(\alpha_{#1-1},\alpha_{#1})}}
\newcommand{\paramijk}[4]{\parambase_{#1#4}^{(#2,#3)}}
\newcommand{\phik}[2]{\phi_{#1#2}^{(\alpha_{#1-1},\alpha_{#1})}}
\newcommand{\phikk}[2]{\phi_{#1#2}}

\newcommand{\fcmidelem}[3]{
    \left[
     \begin{array}{ccc}
       \ffiber{#1}{1,1}(x_{#1}) & \cdots & \ffiber{#1}{1,#3}(x_{#1}) \\
        \vdots & & \vdots \\
       \ffiber{#1}{{#2},1}(x_{#1}) & \cdots & \ffiber{#1}{#2,#3}(x_{#1})
     \end{array}
    \right]}

\newcommand{\tensor}[1]{\bm{\mathcal{#1}}}
\newcommand{\TTr}{\mathbb{T}\mathbb{T}}
\newcommand{\FTr}{\mathbb{F}\mathbb{T}}

\newcommand{\alr}{\varepsilon} % approximate low rank
\newcommand{\crossdelta}{\delta_{\textrm{cross}}} % stopping tolerance for cross approximation
\newcommand{\roundeps}{\epsilon_{\textrm{round}}} % stopping tolerance for cross approximation
\newcommand{\approxeps}{\epsilon_{\textrm{approx}}} % stopping tolerance for cross approximation

\begin{document}

\begin{frontmatter}
  
\begin{abstract}
  We develop new approximation algorithms and data structures for representing and computing with multivariate functions using the functional tensor-train (FT), a continuous extension of the tensor-train (TT) decomposition. The FT represents functions using a tensor-train ansatz by replacing the three-dimensional TT cores with univariate matrix-valued functions. The main contribution of this paper is a framework to compute the FT that employs adaptive approximations of univariate fibers, and that is not tied to any tensorized discretization. The algorithm can be coupled with any univariate linear or nonlinear approximation procedure. We demonstrate that this approach can generate multivariate function approximations that are several orders of magnitude more accurate, for the same cost, than those based on the conventional approach of compressing the coefficient tensor of a tensor-product basis. Our approach is in the spirit of other continuous computation packages such as Chebfun, and yields an algorithm which requires the computation of ``continuous'' matrix factorizations such as the LU and QR decompositions of vector-valued functions. To support these developments, we describe continuous versions of an \revise{approximate maximum-volume} cross approximation algorithm and of a rounding algorithm that re-approximates an FT by one of lower ranks. We demonstrate that our technique improves accuracy and robustness, compared to TT and quantics-TT approaches with fixed parameterizations,  of high-dimensional integration, differentiation, and approximation of functions with local features such as discontinuities and other nonlinearities.
\end{abstract}

% keywords for siam

\end{frontmatter}

\section{Introduction}

Tensor decompositions are a popular tool for mitigating the curse of dimensionality in applications that involve large-scale multiway arrays. Examples are wide-ranging and include compressing the results of scientific simulations~\cite{Austin2016}, neural networks~\cite{Novikov2015}, numerical solution of PDEs~\cite{Khoromskij2010}, stochastic optimal control~\cite{Gorodetsky2015,Gorodetsky2018}, and uncertainty quantification~\cite{Bigoni2016, Gorodetsky2017filter,Eldred2018}. %%

Yet many of the above problems are not naturally posed in the setting of discrete, structured grids. Instead, they are more naturally cast as questions of \textit{function approximation}. Here, low-rank tensor decompositions have been applied by either (i) discretizing and evaluating a multivariate function on a tensor-product grid to form a multiway array~\cite{Khoromskij2011, Dolgov2013, Lee2016}, or (ii) representing the function in a tensor-product basis and decomposing the coefficient tensor~\cite{Mathelin2014, Rai2014}.
But these two methods are not flexible enough for many applications; more general ``continuous'' or functional low-rank representations are needed. In this paper we describe an \emph{adaptive} low-rank function approximation framework that is not limited to particular discretizations or fixed parameterizations.

As an example, consider representing a function with localized features---e.g., a probability density function that concentrates around one or more modes or a function with jump discontinuities; doing so can require specialized (non tensor-product) discretizations to avoid excessive
computation. Indeed, adaptive mesh refinement is an important area of research precisely because efficient and accurate solutions to real-world problems must often  resolve
local structure.
Another limitation, specifically of the discretization approach, is that its utility diminishes outside of the tasks of compressed storage and evaluation.  
To illustrate, consider computing a derivative of a compressed function. If the compressed representation is built from pointwise evaluations on a tensor-product grid, derivative approximations are essentially limited to finite difference rules that can be defined on the same underlying grid.
Yet scalable algorithms for uncertainty quantification or control might first require approximating 
a function, and then operating (e.g., performing differentiation) with the approximation. For instance, in the context of an iterative PDE solver, one might need to apply a Poisson operator with an inhomogeneous coefficient field onto a test function~\cite{mythesis}. In this case, the grid representing the field may differ from the grid representing the test function, and computing inner products or elementwise products of arrays of different sizes is an ill-defined operation. %% 

On the other hand, parameterizing a function  using a prescribed tensor-product basis and then compressing the coefficients is often insufficiently expressive. Modern applications of machine learning suggest that tensor-product bases might not efficiently represent many phenomena; nonlinear parameterizations, e.g., neural networks, can be more successful. In uncertainty quantification and control theory, there has also been a realization that propagating uncertainty, performing inference, or representing feedback policies
for complex models requires a certain degree of nonlinear parameterization, for example, through local approximations~\cite{Conrad2016} or with neural networks~\cite{Bertsekas1995}. %% 

% 
% %
% %
% 
% %
% 

Rather than compressing function evaluations on a fine grid or the coefficients of a prescribed tensor-product basis, we seek a more adaptive and expressive approach, suited to a wide range of computational operations.
The central element of our approach is a continuous analogue of the tensor-train representation, called the \emph{functional tensor-train}, described in~\cite{Oseledets2013}. In that work, a multivariate function is represented as a particular sum of products of univariate functions:
\begin{align}\label{eq:ft}
f(x_1,\ldots,x_d) = \sum_{\alpha_0=1}^{r_0} \sum_{\alpha_1=1}^{r_1} \cdots \sum_{\alpha_d=1}^{r_{d}}\ffiber{1}{\alpha_0,\alpha_1}(x_1)\ffiber{2}{\alpha_1,\alpha_2}(x_2) \cdots \ffiber{d}{\alpha_{d-1},\alpha_{d}}(x_d),
\end{align}
with $r_0 = r_d = 1$. Analytical examples in~\cite{Oseledets2013} demonstrate how certain functions can be represented in this continuous low-rank format. 

In this paper we provide a \textit{computational methodology} for approximating multivariate functions and computing with them in this format. We do not discretize functions on tensor-product grids and we do not \textit{a priori} specify a tensor-product basis for approximation, and as a result our approach overcomes the limitations described above. Furthermore, our approach is adaptive and akin to incorporating adaptive grid refinement within the context of low-rank approximations. 
\revise{While some related approaches~\cite{Beylkin2009,Doostan2013,Gorodetsky2018b} represent low-rank functions without relying on tensor-product discretizations or coefficient tensors, the firs two~\cite{Beylkin2009,Doostan2013} employ the canonical polyadic format and all of them rely on fixed (non-adaptive) parameterizations to learn from scattered data.}
To the best of our knowledge, there are no existing methods for low-rank function approximaton that adapt parameterizations (on a grid or otherwise) in addition to ranks.
Overall, our contributions include:
\begin{enumerate}
\item \textit{locally and globally adaptive} algorithms for low-rank approximation %% 
 of black-box functions;
\item algorithms for rounding multivariate functions already in FT format; and
\item a new \textit{data structure} for representing~\eqref{eq:ft} that is suited to %% 
  both linear and nonlinear parameterizations of the univariate functions.
\end{enumerate}

These contributions are enabled by integrating two lines of research: tensor decompositions and continuous linear algebra. We combine the idea of separation of variables, which underlies tensor decompositions, with the flexibility of linear and multilinear algebraic algorithms for ``computing with functions.'' Computing with functions is an emerging alternative to \textit{discretize-then-solve} methodologies~\cite{Olver2014}, and  has been facilitated by continuous extensions to the common and widely used linear algebra techniques that underpin virtually all numerical simulations. For example, \cite{Battles2004,Platte2010,Townsend2013,Townsend2014,chebfuncode} introduce an extension of MATLAB, called Chebfun, for computing with functions of one, two \cite{Townsend2013}, and three \cite{chebfun3ref} variables. Chebfun implements continuous versions of the LU, QR, and singular value decompositions; enables approximation, optimization, and rootfinding; and provides methods for solving differential equations without discretization. Another example of computing with functions `on-the-fly' is~\cite{Mohlenkamp2013}, where functions are adaptively built, as needed, to solve the Schr{\"o}dinger equation. Our work extends this methodology and portions of Chebfun's capabilities to arbitrary dimensions by marrying them with tensor decompositions.

As demonstrated by Chebfun for almost a decade, the advantages of computing with functions
include increased accuracy and stability, along with more natural modeling and algorithm development: complex operations on functions can be built through relatively simpler operations with great generality. Analogously, our approach provides similar building blocks for multivariate functions. These building blocks lead to new capabilities relative to existing tensor-based methods, including the quantics tensor-train (QTT) decomposition~\cite{Khoromskij2010, Khoromskij2011qtt}; capabilities include progressive \textit{refinement of local features}, \textit{adaptive integration and differentiation} rules that offer significant gains in efficiency and numerical stability, and the ability to apply binary operations to functions that have entirely \textit{different parameterizations.}

The rest of the paper is organized in four sections. In Section~\ref{sec:background}, we detail how tensors appear in function approximation problems and describe the resulting limitations. In Section~\ref{sec:datastruct}, we describe a data structure for representing~\eqref{eq:ft} that stores parameters of the univariate functions \textit{independently}. The structure does not store any tensor, either explicitly or implicitly. We also show how storing a tensor-train decomposition of a coefficient tensor is a specific realization of our framework and comment on the additional flexibility that we add. In Section~\ref{sec:lowrankcompression} we describe continuous analogues of cross-approximation and rounding that are used to decompose a black box function into its low-rank representation. Finally, in Section~\ref{sec:numerics} we present numerical experiments validating our approach that demonstrate significant improvements over previous low-rank functional representations.

\section{Discrete tensor-trains and function approximation}\label{sec:background}
In this section, we provide background on tensor decompositions, review their existing use for representing functions, and describe their limitations. 
Let $\posint$ denote the set of positive integers and $\reals$ the set of reals. Let $d \in \posint$ and $\xspace{} =  \xspace{1} \times \cdots \times \xspace{d}$ be a tensor-product space with $\xspace{k} \subset \reals$ for $k = 1,\ldots,d$. Let $n_{k} \in \posint$ for $k = 1,\ldots,d$ and $N = \prod_{k=1}^dn_k$. A tensor is a $d$-way array with $n_k$ elements along \textit{mode} $k$, and is denoted by uppercase bold calligraphic letters $\tensor{A} \in \reals^{n_1 \times \cdots \times n_d}$. The special case of a 2-way array (matrix) is denoted by an uppercase bold letter, e.g.,  $\mat{A} \in \reals^{m \times n}$, and a 1-way array (vector) is denoted by a lowercase bold letter, e.g., $\bvec{a} \in \reals^{m}$.

\subsection{Tensor-train representation}\label{sec:tt}

A \textit{tensor-train} (TT) representation of a tensor $\tensor{A}$ is defined by a list of $3$-way arrays, $\TTr(\tensor{A}) = \left(\tensor{A}_k \in \reals^{r_{k-1} \times n_k \times r_{k}}\right)_{k=1}^{d}$, with $r_0=r_d=1$ and $r_k \in \posint$ for $k=2,\ldots,d-1$. Each $\tensor{A}_k$ is called the \textit{TT-core}, and the sizes $\left(r_k\right)_{k=0}^d$ are called the \textit{TT-ranks}. Computing with a tensor in TT format requires computing with its cores. For example, an element of a tensor is obtained through multiplication
\begin{equation*}
\tensor{A}[i_1,i_2,\ldots,i_d] = \tensor{A}_1[1,i_1,:]\tensor{A}_2[:,i_2,:]\cdots \tensor{A}_d[:,i_d,1], \quad  1 < i_k < n_k \textrm{ for all } k,
\end{equation*}
and tensor-vector contraction, with $\left( \bvec{w}_k \in \reals^{n_k}\right)_{k=1}^d$, is obtained according to
\begin{equation*}
\tensor{A} \times_1 \bvec{w}_1 \cdots \times_d \bvec{w}_d = \left(\sum_{i_1=1}^{n_k} \tensor{A}_1[1,i_1,:]\bvec{w}[i_1]\right)\left(\sum_{i_2=1}^{n_2} \tensor{A}_2[:,i_2,:]\bvec{w}[i_2]\right)\cdots \left(\sum_{i_d=1}^{n_d}\tensor{A}_d[:,i_d,1]\bvec{w}[i_d]\right).
\end{equation*}
The advantage of this representation of tensors is reduced storage and computational costs that scale \textit{linearly} with dimension, rather than exponentially. If we let $r_k=r$ and $n_k=n$, then the cost of storing the list of cores is $\mathcal{O}(dnr^2)$, and the costs of evaluation and contraction are $\mathcal{O}(dr^2)$ and $\mathcal{O}(dnr^2)$, respectively.\footnote{Contraction of the full tensor $\displaystyle{\sum_{i_1=1}^{n_1}\sum_{i_2=1}^{n_2}\cdots\sum_{i_d=1}^{n_d}\tensor{A}[i_1,i_2,\ldots,i_d]\bvec{w}[i_1]\bvec{w}[i_2]\cdots\bvec{w}[i_d]}$ requires an exponentially growing, $N = \mathcal{O}(n^d)$, number of operations.} The computational complexity becomes driven by the problem's rank.

In practice, tensors  typically have low-rank representations that are \textit{approximate} rather than exact. To make the concept of approximate low-rank representations more concrete, the TT-ranks are related to SVD ranks of various tensor reshapings. Let $\mat{A}^k$ represent a reshaping of the tensor $\tensor{A}$ along the $k$th mode into a matrix: in MATLAB notation we have
\begin{align*}
\mat{A}^k &= \texttt{reshape}\left(\tensor{A},\left[ \prod_{i=1}^{k}n_i, \prod_{i={k+1}}^{d}n_k\right]\right).
\end{align*}
Now consider the SVD of this matrix. According to~\cite{Oseledets2011}% \cite{Thm.\ 2.2}{Oseledets2011}
, if for each unfolding $k$ we have  $\mat{A}^k = \mat{G}^k + \mat{E}^k$ such that the matrices $\mat{G}^k$ have SVD rank equal to $r_k$ and $\lVert \mat{E}^k \rVert = \alr_k$, then a rank $\bvec{r} = (1,r_1,\ldots,r_{d-1},1)$ approximation  $\tensor{A}_{\bvec{r}}$ exists with bounded error
\begin{equation*}
  \lVert \tensor{A} - \tensor{A}_{\bvec{r}} \rVert^2  \leq \sum_{k=1}^{d-1} \alr_k^2.
\end{equation*}
Essentially, this result implies that for problems with good separability, as defined by the matrix SVD, between groups of the first $k$ and last $(d-k)$ variables, an approximate tensor-train representation with low ranks may attain good accuracy.

Finally, we comment on how existing literature converts function approximation problems into tensor decomposition problems. The first method is straightforward: each $\xspace{k}$ is discretized into $n_k$ elements so that $\xspace{}$ is a tensor product grid. The tensor arises from evaluating a multivariate function $f(x_1,\ldots,x_d)$ at each node in this grid. In this case we view $\tensor{A}$ as a mapping from a discrete set to the real numbers, $\tensor{A}:\xspace{} \to \reals$. 

The second method for obtaining a tensor is to represent a function in a tensor product basis, i.e., a basis produced by taking a full tensor product of univariate functions~\cite{Rai2014,Chevreuil2015,Bigoni2016}. In particular, let $$\left\{\phikk{k}{i_k} :\xspace{k} \to \reals: k = 1,\ldots,d, \textrm{ and } i_k = 1,\ldots, n_k\right\}$$ be a set of univariate basis functions. Then a function represented in this basis can be written as
\begin{equation}\label{eq:lowrankbasis}
  f(x_1,x_2,\ldots,x_d) = \sum_{i_1}^{n_1} \sum_{i_2}^{n_2}\cdots \sum_{i_d}^{n_n} \tensor{A}[i_1,i_2,\ldots,i_d] \phikk{1}{i_1}(x_1) \phikk{2}{i_2}(x_2) \cdots \phikk{d}{i_d}(x_d).
\end{equation}
Now the coefficients of this expansion comprise the tensor of interest. This representation arises quite frequently in uncertainty quantification, e.g., polynomial chaos expansions where the univariate basis functions are orthonormal polynomials.\footnote{While the grid discretization approach, i.e., pointwise evaluations of a function on a tensor product grid, could also be thought of as, for instance, a piecewise constant approximation, the latter involves the \emph{additional} assumption of a particular basis or interpolation scheme. Thus we prefer to distinguish the grid discretization approach from the alternative of choosing a basis and storing its coefficients.}

Finally, we point out a simple extension to the tensor-train that is typically more efficient when the $n_k$ are large. The \emph{quantics tensor-train} (QTT)~\cite{Khoromskij2011qtt,Khoromskij2011} decomposition increases the dimension of a tensor by reshaping the standard modes into $2 \times 2 \times \cdots$ modes,  $3 \times 3 \times \cdots$ modes, or larger-sized modes. Typically such a reshaping allows a further reduction in storage complexity so that it scales with $\log(n)$ rather than $n$.

\subsection{Limitations of existing tensor-train approaches}
In this section, we describe some limitations of the TT representation for function approximation, which follow from its discrete nature; these limitations will motivate continuous extensions. Among our broader goals is to use low rank approximations to mitigate the curse of dimensionality in uncertainty quantification and control.
Both application areas share procedures that employ simple operations---addition, multiplication, integration, differentiation, and computing inner products---to construct more complex algorithms. 
However, the direct application of existing low-rank tensor approximations within such algorithms is faced with the following challenges:
\begin{enumerate}
\item Employing low-rank tensors in new contexts 
 requires choosing specific discretizations or parameterizations (e.g., what grid to use, what quadrature rule to incorporate); the fact that these choices are not automatic limits generality and adaptivity.
\item Operating on functions with different discretizations requires ad hoc %%
interpolation onto shared grids or a shared parameterization, before executing an operation. %% 
\end{enumerate}
Specific limitations associated with these challenges are summarized in Table~\ref{tab:tensorlimitations}. We illustrate these limitations with three examples:
\begin{table}
\centering
    \begin{tabular}{ccc}
      \toprule
      Desired Usage & Grid discretization & Tensor-of-coefficients \\
 \midrule
      Represent local features & no adaptivity & only linear parameterizations\\ & & and no adaptivity \\ 
      Extract gradients & only finite differences & \revise{\textbf{no limitation}} \\[0.2em]
      Perform integration & only tensor-product quadrature & \revise{\textbf{no limitation}} \\[0.2em]
      Multilinear algebra with & only identical grids & only identical\\ several functions & & parameterizations \\
\bottomrule
    \end{tabular}
    \caption{Limitations of grid discretization and tensor-of-coefficients formulations for computing with functions} 
    \label{tab:tensorlimitations}
    \vspace{-15pt}
\end{table}

\textit{(1) Approximating functions with discontinuities or local features.} Multivariate functions with localized nonlinearities or discontinuities arise in, among other settings, stochastic optimal control problems with discontinuous cost functions~\cite{Gorodetsky2018}. In these cases, specifying basis functions \textit{a priori} and approximating the resulting tensor of coefficients can lead to numerical issues: for instance, if the basis functions do not respect the discontinuity, the Gibbs phenomenon can ruin approximation quality. On the other hand, a fixed wavelet basis can potentially under-resolve or over-resolve the discontinuity. To the best of our knowledge, adaptive schemes to address these issues in the setting of high-dimensional low-rank function approximation %%
do not yet exist.

\textit{(2) Computing gradients and integrals of multivariate functions.} 
Gradients and integrals appear within a myriad of algorithms, and efficient ways to compute them in high dimensions are essential. For example, value and policy iteration in dynamic programming~\cite{Bertsekas2007} involve updating a multivariate function through fixed point iterations, and these updates require the gradient of a multivariate function~\cite{Fleming2006}. If we simply represent a multivariate function on a discretized grid, obtaining the gradient is an undefined operation. Numerical finite difference schemes to approximate gradients face the problem of choosing a good discretization. Dealing with such issues by, for example, discretizing each input coordinate %%
with an extremely large number of points % (e.g., $2^{25}$) 
and performing a QTT decomposition does not solve this problem, as finite-difference roundoff errors will persist. %% 

\textit{(3) Computing with several functions}.
Once a low-rank approximation of a multivariate function is constructed, there are enormous computational advantages to \textit{computing} with the function in this low-rank format. For example, consider multiplying or adding two functions. In the discrete setting of a tensor obtained from pointwise evaluations of a function on a tensor product grid, these operations are undefined unless both functions have identical underlying discretizations. Otherwise, problem-specific interpolation techniques are needed in order to interpolate the functions onto the same grid. These interpolation procedures are potentially inefficient and can result in unnecessary decompression of the underlying function. Similar procedures would be needed to compute identical parameterizations for the case of a tensor of coefficients. For example, suppose that two functions are parameterized by piecewise polynomials with different choices of knots: existing algorithms would require a projection or interpolation of both functions onto the same parameterization.

The limitations above arise from a \textit{data structure} problem. When we represent a multiway array of function values in TT or QTT format, the connection between these values and the function itself is lost. In other words, one is left to create algorithms that work with discrete objects without any knowledge of where these objects came from. One must then inject additional knowledge in a post-processing stage to perform analysis. Consider, for example, the problem of integrating a multivariate function. Given function values on a tensor of input points, one could perform Monte Carlo integration, apply a Newton--Cotes rule, perform rational interpolation and then integrate, etc. Each of these choices involves assumptions about the properties of the underlying function that are not available in the discrete tensor representation itself.
Ideally, a data structure that combines \textit{function values or parameters} and methods that \textit{interpret} them, through both the compression and post-processing steps, are needed. 

While each of the limitations and examples provided above can potentially be addressed independently, we seek an algorithmic framework that is general enough to overcome all of them. Such a framework requires both a new data structure that more closely follows the mathematical model~\eqref{eq:ft}, and an abstraction of discrete algorithms to the continuous case. Fortunately, Chebfun~\cite{Platte2010} and related efforts have formulated exactly such a computational paradigm. In fact, they have demonstrated the great success of such approaches for over ten years, for problems ranging from function approximation to solving PDEs. The rest of this paper describes how to apply this paradigm to low-rank functions of arbitrary dimension.

\section{Continuous analogue of the tensor-train}\label{sec:datastruct}

Our proposed data structure 
to represent the functional tensor-train~\eqref{eq:ft} stores each univariate function $\ffiber{k}{\alpha_{k-1},\alpha_k}$ \textit{independently} and groups the functions associated with each input coordinate $k \in \{1, \ldots, d \}$ into matrix-valued functions $\mvf{F}_k:\xspace{k} \to \reals^{r_{k-1}\times r_{k}}$, called cores. Note that no three-way arrays %%
 are required. The evaluation of a function in this format, using these continuous cores, involves multiplying matrix-valued functions:
\begin{equation}\label{eq:ft_with_cores}
f(x_1,x_2,\ldots,x_d) = \fcore{1}\fcore{2} \cdots \fcore{d}.
\end{equation}
In this section, we describe a new parameterization of this format and compare it to existing work.

Working with this mathematical object on a computer requires a finite parameterization. In our case, we parameterize~\eqref{eq:ft} with a hierarchical structure. Beginning at the leaf level, %%
 each univariate function $\ffiber{k}{\alpha_{k-1},\alpha_k}$ 
is parameterized by $\fn{k} \in \posint$ parameters $\left(\paramk{k}{j}\right)_{j=1}^{\fn{k}}$. Both linear parameterizations, e.g., a basis expansion
\begin{equation}\label{eq:linbasis}
\ffiber{k}{\alpha_{k-1},\alpha_{k}}(x_k) = \sum_{j=1}^{\fn{k}}\paramk{k}{j}\phik{k}{j}(x_k),
\end{equation}
and nonlinear parameterizations are allowed. As an example  of the latter, consider a  piecewise constant approximation,
\begin{equation}\label{eq:piecewise}
\ffiber{k}{\alpha_{k-1},\alpha_{k}}(x_k) = \left\{
  \begin{array}{cl} 
    \paramk{k}{(j+n/2)} & \textrm{ if }  \paramk{k}{j} \leq x < \paramk{k}{j+1} \\
    0 & \textrm{otherwise}
  \end{array}
  \right. ;
\end{equation} %
where the parameters are here a concatenation of the knot locations and the function value between any pair of knots. In practice, we employ even more adaptive nonlinear parameterizations---refining piecewise polynomial approximations, and the value of $\fn{k}$, independently for each univariate function. Similarly, the number of basis functions in the linear parameterization \eqref{eq:linbasis} can even vary among the univariate functions for a given $x_k$. In Section~\ref{sec:gaussbump} we show a numerical example where both linear and nonlinear parameterizations can be used \textit{simultaneously} for an input coordinate $x_k$.

Note that these parameterizations cannot be incorporated within the tensor-of-coefficients framework. As a result, our approach is more expressive precisely because it decouples the notion of ``tensor products of functions'' from the notion of ``tensors of parameters.''
We also note that while any function \textit{can} be represented using a complete tensor-product basis, this is seldom the \textit{most efficient} choice. Such representations are often convenient for proving theoretical results, but they may be too computationally expensive in practice.

The second level of the hierarchy groups $\ffiber{k}{\alpha_{k-1},\alpha_{k}}:\xspace{k} \to \reals$ into a matrix-valued function $\mvf{F}_k:\xspace{k} \to \reals^{r_{k-1} \times r_{k}}$, i.e.,
\begin{equation}\label{eq:ftcore}
\mvf{F}_k(x_k) = \fcmidelem{k}{r_{k-1}}{r_{k}}.
\end{equation}

The third level, analogous to the discrete setting, represents the continuous tensor-train with a list of matrix-valued functions, $\FTr(f) = \left(\mvf{F}_k:\xspace{k} \to \reals^{r_{k-1} \times r_{k}}\right)_{k=1}^{d}$, with $r_0=r_d=1$ and $r_k \in \posint$ for $k=2,\ldots,d-1$, as described by~\eqref{eq:ft_with_cores}.

\begin{figure}
  \includegraphics[width=\textwidth]{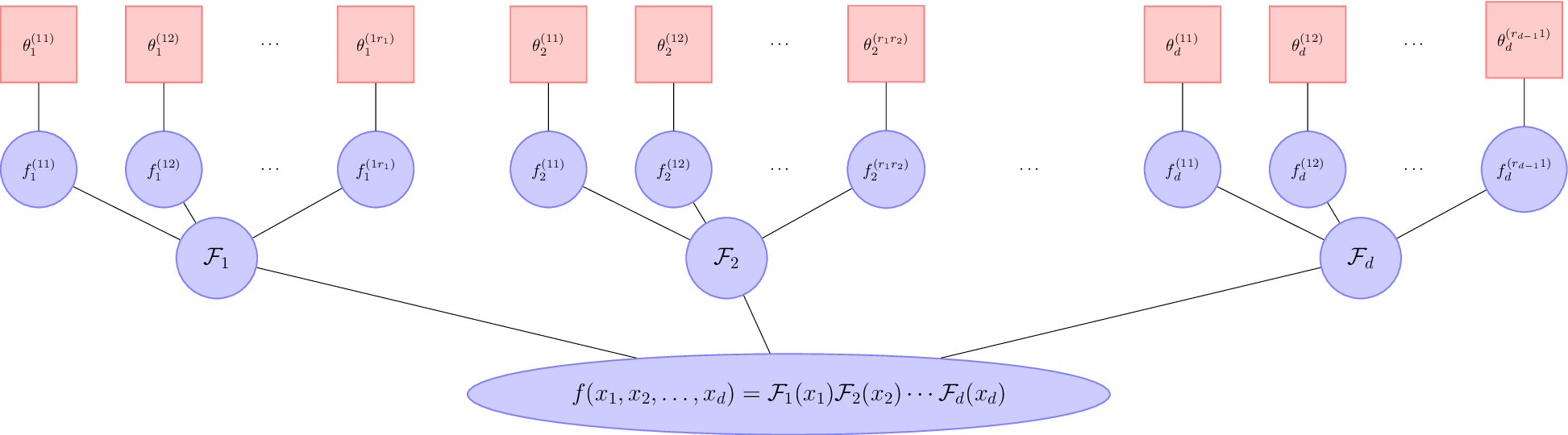}
  \caption{Hierarchical data structure for storing low-rank functions. Parameters are represented by red squares and functions are represented in blue. In contrast to previous approaches~\cite{Chevreuil2015,Bigoni2016}, the parameters and parameterized form that represent each univariate function are stored separately, i.e., no tensor-product structure is imposed upon the \textbf{parameters} of the constituent univariate functions. Instead tensor-product structure is imposed upon the \textbf{univariate functions} themselves.}
  \label{fig:fthierarchy}
\end{figure}

A graph %% 
representation of this hierarchy is shown in Figure~\ref{fig:fthierarchy}. There are  \textit{no tensor data structures} in this representation; there is no underlying tensorized grid or tensor of coefficients, either explicitly or in a decomposed form.  Next we demonstrate conditions under which our representation is equivalent \revise{to~\eqref{eq:lowrankbasis}} and when it is advantageous.

\subsection{Comparison with a tensor of coefficients}

We now compare the above approach with the typical alternative of representing multivariate functions via a tensor-train decomposition of coefficients. In particular, we demonstrate conditions under which the low-rank representation~\eqref{eq:ft} and a low-rank representation of a coefficient tensor in~\eqref{eq:lowrankbasis} lead to equivalent approximations. This comparison will serve to highlight the advantages of directly computing with~\eqref{eq:ft} instead of~\eqref{eq:lowrankbasis}.

\subsubsection{Conditions of equivalence}

Equivalence between the two models arises when the following two conditions are met:
\begin{enumerate}
  \item Each $\ffiber{k}{\alpha_{k-1},\alpha_{k}}$ is parameterized linearly with a fixed number of basis functions in~\eqref{eq:linbasis}.
  \item An \textit{identical} basis is used for all $\ffiber{k}{\alpha_{k-1},\alpha_{k}}$ within a dimension $k$, i.e., $\fn{k} = \fnn{k}$ and $\phik{k}{j} = \phikk{k}{j}$ for all $\alpha_{k-1} = 1, \ldots, r_{k-1}$, $\alpha_{k} = 1, \ldots, r_{k}$, and $j = 1 \ldots, \fnn{k}.$
\end{enumerate}

To demonstrate this fact, consider a tensor $\tensor{A}_k \in \reals^{r_{k-1} \times \fnn{k} \times r_{k}}$ such that
\begin{equation*}
\tensor{A}_k[:,j,:] =
\left[
     \begin{array}{ccc}
        \paramijk{k}{1}{1}{j} & \cdots & \paramijk{k}{1}{r_k}{j} \\
        \vdots & & \vdots \\
       \paramijk{k}{r_{k-1}}{1}{j} & \cdots & \paramijk{k}{r_{k-1}}{r_k}{j} \\
     \end{array}
    \right],
\end{equation*}
for $j = 1, \ldots, \fnn{k}$. This tensor is the $k$th core of a TT decomposition of the coefficient tensor in~\eqref{eq:lowrankbasis}. The relationship between the continuous core and the discrete TT core is now
\begin{equation*}
  \mvf{F}_k(x_k) = \sum_{j=1}^{\fnn{k}}\tensor{A}_k[:,j,:]\phikk{k}{j}(x_k),
\end{equation*}
where the same basis function is used for every slice of $\tensor{A}_k$. In other words, one contracts the TT core with the basis functions along the second mode. This connection yields the equivalence:
\begin{align*}
f(x_1,\ldots,x_2) = \mvf{F}_1(x_1)\cdots \mvf{F}_d(x_d) 
                  = \sum_{i_1=1}^{\fnn{1}}\cdots\sum_{i_d=1}^{\fnn{d}} \tensor{A}_1[:,i_1,:] \cdots \tensor{A}_d[:,i_d,:] \phikk{1}{i_1}(x_1) \cdots \phikk{d}{i_d}(x_d)
\end{align*}

\subsubsection{FT as a sparse storage structure}

One advantage of our approach is that, under the equivalence conditions described above, it can be viewed as a \textit{sparse} storage scheme for the TT cores $\tensor{A}_k$. We demonstrate this property on two examples.

The first example involves representing additive functions. Additive functions are extremely common within high-dimensional modeling~\cite{Hastie1990,Lie2008,Meier2009}. As such, it is useful to be able to represent these functions in low-rank format. An additive function has TT-ranks $r_1= \cdots = r_{d-1} = 2$, which can be seen from the following decomposition:
\begin{equation*}
  f(x_1,x_2,\ldots,x_d) = f_1(x_1) + f_2(x_2) + \cdots + f_d(x_d) 
                       = \left[f_1(x_1) \  1\right] \left[ 
                         \begin{array}{cc}
                           1 & 0 \\
                           f_2(x_2) & 1 
                         \end{array} \right] \cdots 
                                      \left[ \begin{array}{c}
                                        1 \\ 
                                        f_d(x_d) 
                                      \end{array}\right].
\end{equation*}
Assume that $\fnn{k} = n$ terms are required to parameterize each $f_k$ in a basis expansion. Then the TT cores of the coefficient tensor store $4n + 4(d-2)n$ floating point numbers. The constants $1$ and $0$ are essentially represented with $n$ parameters, when in reality they can be represented with a single parameter. Standard techniques in the literature do not exploit this structure of the TT cores. If we instead let $\fn{k}$ vary for each function in each core (violating the second equivalence condition), only $2(n+1) + (d-2)(n+3)$ floating point numbers would require storage. In other words, in the limit of large $d$, four times less memory is required for our low-rank format than for a TT decomposition of coefficients.

A more impressive example is a quadratic function. A quadratic function has $d^2$ parameters, and therefore one would expect any low rank storage format only to require storing $\mathcal{O}(d^2)$ parameters. However, this storage requirement can only be achieved by a sparse storage scheme for the TT cores in~\eqref{eq:linbasis}, which is naturally supported by the  format that we propose. To be more clear, one parameterization of a quadratic function $$f(x_1,x_2,\ldots,x_d) = x^T\mat{A}x,\quad  x = [x_1 \ x_2 \ \cdots \ x_d]$$ in low-rank format has continuous cores with the following structure:\footnote{For clarity of presentation we have written these cores as being $d \times d$, i.e., a rank $d$ representation. A more compact representation is possible with cores of varying sizes and the maximum core size being $d/2 \times d/2$. The asymptotic results that we present would not change in this case, however.}
\begin{align*}
\mvf{F}_1(x_1) =  \left[ 
\begin{array}{c@{\hspace{2em}}cccc}
a_{1,1} x_1^2 & a_{1,2}x_1  & \ldots & a_{1,d} x_1  & 1
\end{array}
\right], \qquad 
\mvf{F}_d(x_d) = \left[
\begin{array}{ccc}
1 &
x_d &
a_{d,d}x_d^2 
\end{array}
\right]
\end{align*}

\begin{align*}
\mvf{F}_k(x_k) =  \left[
\begin{array}{c@{\hspace{0.5em}}c@{\hspace{0.5em}}c@{\hspace{0.5em}}c@{\hspace{0.5em}}c@{\hspace{0.5em}}c@{\hspace{0.5em}}c@{\hspace{0.5em}}}
1  	     & 0      & 0      & 0      & \cdots & 0      & 0        \\ 
x_k      & 0      & 0      & 0      & \cdots & 0      & 0       \\
0 		 & 1      & 0      & 0      & \cdots      & 0      & 0     \\% [2.5em]
0 		 & 0      & 1      & 0      & \cdots & 0      & 0       \\ % [1em]
\vdots 	 &        & \ddots & \ddots &        &        & \vdots    \\ % [1.5em]
\vdots 	 &        &        & \ddots & \ddots &        & \vdots      \\ % [2.5em]
\vdots 	 & \vdots &        &        &        & \ddots & \ddots \\ %[1em]
0 		 &  0     & 0      & 0      & \cdots & 1      & 0       \\
a_{k,k}x_k^2  & a_{k,k+1}x_k  & a_{k,k+2}x_k & \cdots & a_{k,d-1} x_k & a_{k,d} x_k & 1
\end{array}
\right] .
\end{align*}

The TT representation of each core requires storing $d \times 3 \times d$ floating point values, denoting the coefficients of constant, linear, and quadratic basis functions. The data structure that we propose instead stores each element of $\mvf{F}_k(x_k)$ separately, and there are $\mathcal{O}(d)$ nonzero elements. Thus the total storage of our proposed structure is $\mathcal{O}(d^2)$, versus $\mathcal{O}(d^3)$ for the TT version. 

Similar examples can be constructed for functions of higher polynomial degree. The main takeaway from these examples is that the data structure that we propose is flexible enough to take advantage of core sparsity. Thus, even when it is appropriate to model a function using~\eqref{eq:linbasis}, we can achieve better compression using~\eqref{eq:ft}. Secondly, we note that while it is true that one could potentially use a sparse storage mechanism for the TT cores to achieve similar benefits, our proposed framework naturally supports this structure with no additional modifications.

\subsubsection{Nonlinear parameterizations for more accurate and efficient compression}\label{sec:nonlinparam}

This section provides an example of a function that, while representable in a low-rank format, cannot be readily described by a decomposition of a coefficient tensor. In particular, we demonstrate the additional expressivity afforded by modeling functions using~\eqref{eq:ft}. % instead of~\eqref{eq:linbasis}.

Consider the rank-2 function
\begin{equation*}
f(x_1,x_2) = \left\{ 
\begin{array}{cc}
0 & \textrm {if }  \frac{1}{4} \leq x_1 \leq \frac{3}{4} \textrm{ and } \frac{1}{4} \leq x_2 \leq \frac{3}{4} \\
\sin(x_1)\cos(x_2) & \textrm{ otherwise}
\end{array}
\right.,  \quad (x_1,x_2) \in [0,1]^2.
\end{equation*}
The rank-2 representation of this function is
\begin{equation}
f(x_1,x_2) = \sin(x_1)\cos(x_2) - \chi_{1/4\leq x_1 \leq 3/4}\sin(x_1)\chi_{1/4\leq x_2 \leq 3/4}\cos(x_2),
\end{equation}
where $\chi_A$ denotes an indicator function on the set $A$.

This function is a model problem for situations when there is a ``hole''  in an otherwise tensor-product space. For example, this might be a rectangle around which a fluid is flowing, an obstacle around which a robot must navigate in a motion planning problem, or a failure region for some black-box model. Sometimes, the location of this region is not known \textit{a priori}, but discovered through interacting with the function. In either case, we would like our function to be zero for inputs that are inside of the ``hole.''

When the boundary is exactly known, this case can easily be handled using the tensor-of-coefficients formulation~\eqref{eq:linbasis}. First one defines a set of basis functions for the ``good'' region; in this case suppose we use $\phikk{1}{1} = \sin(x_1)$. Then for each basis function, one simply adds another basis function according to the rule $\phikk{1}{(2j)} = \chi_{A}\phikk{1}{j}$; in this case we would use $\phikk{1}{2} = \chi_{1/4\leq x_1 \leq 3/4}\sin(x_1).$ Thus, one simply doubles the number of basis functions in each dimension.

In practice, however, the exact location of such failure/zero regions may not be known. Instead, it must be found through the process of evaluating $f$. One way to model such regions is to adapt the knot locations of a piecewise polynomial; see~\eqref{eq:piecewise}. This situation immediately becomes one where the parameters of each representation are no longer coefficients of a tensor product basis. Section~\ref{sec:gaussbump} will provide a numerical example of such a procedure. %

\subsection{Errors introduced through parameterization}\label{sec:theory}

In practice, approximation of each univariate function $\ffiber{k}{ij}$---via a finite-dimensional parameterization---introduces errors in the overall approximation of $f$. 
In our implementation, these errors are incurred during the cross approximation procedure (Section~\ref{sec:lowrankcompression}), when approximating selected univariate fibers of $f$.
The following result relates the
approximation error of each univariate function in~\eqref{eq:ft} to the approximation error in $f$.

\newcommand{\sumconst}{r^{d-1}}
\begin{theorem}[Parameterization error]
 Assume that each of the univariate functions in~\eqref{eq:ft} is bounded according to $\left| \ffiber{k}{ij}(x_k)\right| \leq C$, \revise{ for some $1 \leq C < \infty$.}
  Let $\fafiber{k}{ij}(x_k)$ be univariate functions such that $\lvert \ffiber{k}{ij}(x_k) - \fafiber{k}{ij}(x_k) \rvert \leq \epsilon C$ for $\epsilon d < 1/e$. 
Then the difference between the multivariate approximation
\begin{equation}\label{eq:ftapprox}
\hat{f}(x_1,\ldots,x_d) = \sum_{\alpha_0=1}^{r_0} \sum_{\alpha_1=1}^{r_1} \cdots \sum_{\alpha_d=1}^{r_{d}}\fafiber{1}{\alpha_0,\alpha_1}(x_1)\fafiber{2}{\alpha_1,\alpha_2}(x_2) \cdots \fafiber{d}{\alpha_{d-1},\alpha_{d}}(x_d)
\end{equation}
and~\eqref{eq:ft} is bounded by
\begin{equation}
  \left | f(x) - \hat{f}(x) \right |  \leq  e d^2 C^d \sumconst \epsilon.
\end{equation}
\end{theorem}

\begin{proof}
   Let $\efiber{k}{\alpha_{k-1}\alpha_{k}}(x) = \ffiber{k}{\alpha_{k-1}\alpha_{k}}(x) - \fafiber{k}{\alpha_{k-1}\alpha_{k}}(x)$. First, substitute the definition of each $\fafiber{k}{ij}$ into~\eqref{eq:ftapprox}:
\begin{align*}
\hat{f}(x) &= \sum_{\bvec{1} \leq \bvec{\alpha} \leq \bvec{r}}\left(\ffiber{1}{\alpha_0\alpha_1} - \efiber{1}{\alpha_0\alpha_1}\right)\cdots \left(\ffiber{d}{\alpha_{d-1}\alpha_{d}} - \efiber{d}{\alpha_{d-1}\alpha_d}\right) \\ 
  &= \sum_{\bvec{1} \leq \bvec{\alpha} \leq \bvec{r}}\left[\sum_{\bvec{k} \leq \bvec{1}_d} \left[\ffiber{1}{\alpha_0\alpha_1}\right]^{k_1}\left[-\efiber{1}{\alpha_0\alpha_1}\right]^{1-k_1} \cdots \left[\ffiber{d}{\alpha_{d-1}\alpha_{d}}\right]^{k_d}\left[-\efiber{d}{\alpha_{d-1}\alpha_{d}}\right]^{1-k_d}\right],
\end{align*}where the second equality follows from the binomial theorem.\footnote{Here we have used multi-index notation to simplify the expression. In multi-index notation we assume $\bvec{\alpha}=\left(\alpha_0,\ldots,\alpha_d\right)$, $\bvec{r} = \left(r_0,\ldots,r_d\right)$, $\bvec{k} = \left(k_0,\ldots,k_d\right)$ and $\bvec{1}_d = \left(1,\ldots,1\right)$. The notation $\bvec{\alpha} \leq \bvec{r}$ implies $\alpha_i \leq r_i$ for every  $i = 0, \ldots, d$.}
Computing the difference between this expression and $f$, the term with $\bvec{k} = \bvec{1}_d$ drops out, leaving the rest of the expression
\begin{align*}
f(x) - \hat{f}(x) 
 = -\sum_{\bvec{1} \leq \bvec{\alpha} \leq \bvec{r}} \left[\sum_{\bvec{k} \leq \bvec{1}_d,\bvec{k} \neq \bvec{1}_d} \left[\ffiber{1}{\alpha_0\alpha_1}\right]^{k_1}\left[-\efiber{1}{\alpha_0\alpha_1}\right]^{1-k_1} \cdots \left[\ffiber{d}{\alpha_{d-1}\alpha_{d}}\right]^{k_d}\left[-\efiber{d}{\alpha_{d-1}\alpha_{d}}\right]^{1-k_d}\right].
\end{align*}
Now using the bound $C$ on the univariate functions and $\epsilon C$ on the $\efiber{i}{\alpha_{i-1}\alpha_{i}}$ results in
\begin{align*}
\lvert f(x) - \hat{f}(x) \rvert
&\leq \left \lvert \sum_{\bvec{1} \leq \bvec{\alpha} \leq \bvec{r}} \left[\sum_{\bvec{k} \leq \bvec{1}_d,\bvec{k} \neq \bvec{1}_d} C^{k_1} \left[\epsilon C\right]^{1-k_1} \cdots C^{k_d} \left[\epsilon C\right]^{1-k_d} \right] \right \rvert  \\
&=\revise{\sum_{\bvec{1} \leq \bvec{\alpha} \leq \bvec{r}} \sum_{\bvec{k} \leq \bvec{1}_d,\bvec{k} \neq \bvec{1}_d} \left[\prod_{l=1}^d\left(C^{k_l}C^{1-k_l}\right)\right] \epsilon^{d-\sum_{l}k_l}} \\
&= \revise{\sum_{\bvec{1} \leq \bvec{\alpha} \leq \bvec{r}} C^d \sum_{\bvec{k} \leq \bvec{1}_d,\bvec{k} \neq \bvec{1}_d} \epsilon^{d-\sum_{l}k_l},}
\end{align*}
where in the second equality we combined exponents with common bases. Using  $\sum_{\bvec{1} \leq \bvec{\alpha} \leq \bvec{r}} C^d \leq C^d\sumconst$, \revise{where $r = \max_{1\leq i \leq d-1} r_i$}, we obtain
\begin{align*}
  \lvert f(x) - \hat{f}(x) \rvert \leq C^d \sumconst \sum_{\bvec{k} \leq \bvec{1}_d,\bvec{k} \neq \bvec{1}_d}  \epsilon^{d-\sum_{l}k_l} .
\end{align*}
The $k_l$ terms only enter through a summation. Thus we can simplify the sum over $\bvec{k}$ by summing over possible values for $\sum_{l=1}^d k_{l}$. This sum ranges from zero to $d-1$ and each particular sum is obtained when a certain number of $k_{l}$ are non-zero. Thus we can convert this summation to a sum over $d$ choose $k$ combinations for $k = 0, \ldots, d-1$ to obtain
\begin{align*}
  \lvert f(x) - \hat{f}(x) \rvert  &\leq C^d \sumconst \sum_{l=0}^{d-1}\binom{d}{l} \epsilon^{d-l} = C^d \sumconst \sum_{l=0}^{d-1}\binom{d}{d-l} \epsilon^{d-l}  \leq C^d \sumconst  \sum_{l=0}^{d-1} \left(ed \epsilon\right)^{d-l} \\
  &\leq C^d \sumconst  \sum_{k=1}^{d} \left(ed \epsilon\right)^{k} \leq  d C^d \sumconst ed \epsilon \leq e d^2C^dr^{d-1}  \epsilon ,
\end{align*}
where in the second inequality we have used an upper bound for $\binom{d}{d-l}$ and in the fourth inequality we used the fact that $ d \epsilon < 1/e$.
\end{proof}

This bound formalizes what may be intuitively obvious: the error in a single univariate function is multiplied against all combinations of all the univariate functions of the other input coordinates. Therefore, the impact of any error in a single univariate function can grow exponentially with dimension. At the same time,  $\lvert f \rvert  \leq C^dr^{d-1}$ is the only bound on $f$ itself that we can obtain without further assumptions on the interactions among the univariate functions $\ffiber{k}{\alpha_{k-1}\alpha_k}.$ This fact can be simply seen by setting $\efiber{k}{\alpha_{k-1}\alpha_k} = \ffiber{k}{\alpha_{k-1}\alpha_k}$ in the proof above. Our result can thus be interpreted, perhaps more naturally, as
\begin{equation}
  \left | f(x) - \hat{f}(x) \right |  \leq e d^2  \epsilon \max_{x} \left[f(x)\right] .
\end{equation}
This bound suggests that, to maintain a fixed bound on the relative $L^{\infty}$
error, a conservative approximation scheme should decrease the error threshold $\epsilon$ for the univariate functions comprising~\eqref{eq:ft} quadratically with dimension. In practice, though, we usually  see good approximation behavior with respect to dimension even with fixed thresholds.
Additional assumptions on the interactions between each of the univariate functions can potentially lead to less restrictive requirements for $\epsilon$. %% \AGnote{potentiall add a new assumption}

\section{Constructing low-rank approximations}\label{sec:lowrankcompression}
Next 
we describe how to construct FT approximations with continuous analogues of cross approximation and rounding algorithms~\cite{Oseledets2010}. There are four important differences between the discrete and continuous versions: 
\begin{itemize}
  \item The continuous algorithms interpret tensor fibers as univariate functions obtained by fixing all but one input coordinate;
  \item Each tensor fiber/univariate function is approximated \textit{adaptively} to a particular tolerance as it is encountered, introducing a new source of error;
  \item Optimization \revise{within an approximate maximum volume procedure} for selecting fibers is performed over a \textit{continuous} rather than \textit{discrete} variable; and 
  \item Rounding algorithms use \textit{continuous} QR decompositions, implicitly defining a continuous, rather than a discrete, inner product.
\end{itemize}
Note that the algorithms we present are also applicable to low-rank approximation techniques such as DMRG-cross~\cite{Savostyanov2011}, in which case, rather than tensor fibers and univariate functions, we would have tensor slices and bivariate functions. In other words, our essential contributions lie in placing these existing algorithms in the context of Chebfun-style continuous computation.
Furthermore, if interpreted in the tensor-of-coefficients context, they are the first to provide an \textit{adapt-to-tolerance} procedure for both ranks \emph{and} fibers. 

\subsection{Cross approximation}

In this section, we describe a continuous analogue of a cross approximation algorithm for compressing multiway arrays. Cross approximation is a dimension-sweeping algorithm that seeks to interpolate a tensor using a basis consisting of finite set of its fibers. In the continuous context, we are approximating a multivariate \textit{function} using certain univariate functions formed by fixing all but a single input variable.

To illustrate how continuous cross approximation differs from discrete cross approximation, it is sufficient to consider a bivariate function approximation problem. Extensions to higher dimensions are made using the same types of modifications we describe here to a multivariate cross approximation algorithm, for example \cite[Algorithm 1]{Savostyanov2011}. In two dimensions, cross approximation represents a low-rank function by its CUR decomposition~\cite{Goreinov1997,Mahoney2009,Boutsidis2014}:
\begin{align*}
  f(x,y) &= \mvf{C}(x) \mat{F}^{\dagger} \mvf{R}(y)\\
  \mvf{C}(x) &= \left[f\left(x,y^{(1)}\right)  \ f\left(x,y^{(2)}\right)  \ \cdots \ f\left(x,y^{(r)}\right)\right] \quad   \mvf{R}(y) = \left[ \begin{array}{c}
                         f\left(x^{(1)},y\right) \\
                         f\left(x^{(2)},y\right) \\
                         \vdots \\
                         f\left(x^{(r)},y\right)
                       \end{array}
                \right]       \\
  \mat{F}[i,j] &= f\left(x^{(i)},y^{(j)}\right) \quad i,j = 1,\ldots,r,
\end{align*}
where $\mat{F}^{\dagger}$ indicates the pseudoinverse of $\mat{F}$.
The difference between continuous and discrete CUR decompositions is that columns and rows define a ``quasimatrix'' (vector-valued function)\footnote{Called a quasimatrix because it corresponds to a matrix of infinite rows and $n$ columns; see, e.g.,~\cite{Battles2004,Townsend2014}.} rather than a matrix. In this bivariate case, the fibers are rows and columns. In the more general multivariate case there are fibers corresponding to each dimension of the tensor. As a result, instead of the three-way arrays that form the discrete TT cores, we obtain matrix-valued functions~\eqref{eq:ft_with_cores}.

In the TT (and DMRG-)cross approximation algorithms, fibers associated with each input coordinate are chosen sequentially to \revise{approximately} maximize the volume of the skeleton matrix $\mat{F}$ formed by their intersection. We refer to~\cite{Goreinov1997maxvol,Goreinov2001,Savostyanov2014} for results regarding the quasioptimality of choosing an actual maximum volume submatrix.  This algorithm can be abstracted to the continuous case by generalizing the concept of tensor fibers to encompass the case of univariate functions obtained by fixing all but one input variable. With this abstraction, continuous cross approximation is almost identical (algebraically) to the discrete case. 

Using the bivariate setting for simplicity, we now point out aspects of the algorithm that differ between the discrete and continuous cases. As mentioned above, these algorithmic changes are directly applicable to higher dimensional contexts as well. Using the notation of matrix-valued functions, pseudocode for our continuous cross approximation algorithm is provided in Algorithm~\ref{alg:crossapprox}. The first difference can be seen in Lines~\ref{alg:crossapprox:line:formC}, \ref{alg:crossapprox:line:formR}, and \ref{alg:crossapprox:line:formC2}, where column and row quasimatrices are formed. In the context of multivariate functions, we do not have access to a discrete vector representing a row. Instead we adaptively approximate each column and row via an \textit{online} procedure. In other words, we only generate an approximation of a tensor fiber when the cross approximation algorithm dictates that we need that particular fiber.

The advantage in expressivity of our approach arises during this approximation step. In particular, we adapt the function at the (local) \textit{fiber level} rather than the (global) \textit{dimension} level. Therefore local features occuring along one fiber will not increase the computational effort to approximate another. In Section~\ref{sec:fiberadapt}, we will provide an example wherein global polynomials are used to represent regions where a function is smooth and adaptive piecewise polynomials are used in a region where a function has discontinuities. 

The second difference involves Lines~\ref{alg:ca:line:qr1} and~\ref{alg:ca:line:qr2}. In these lines a QR decomposition of the rows and columns is performed. This step is used to increase the stability of the algorithm; see~\cite{Oseledets2010}. Because we are operating with quasimatrices, we use a \textit{continuous} QR decomposition~\cite{Townsend2013,Townsend2014,mythesis}. Essentially, the discrete and continuous QR decompositions differ with respect to the inner product used to define orthonormality. In the context of a tensor formed by discretizing a function, the discrete inner product is typically only an \textit{approximation} of the continuous inner product. Here, we see that the continuous computation paradigm of Chebfun enables maintaining a notion of orthonormality that is consistent with the original function space, as compared to a discretized TT/QTT approach.

The final difference involves the \texttt{dominant} routine. This routine chooses the set of fibers to be used within a dimension during the next sweep across dimensions. We have developed a continuous analogue of \texttt{dominant} that operates on matrix-valued functions and that is applicable to both the bivariate and multivariate cases. This algorithm is described in the next section.

\begin{algorithm}[ht!]
\renewcommand{\algorithmicrequire}{\textbf{Input:}}
\renewcommand{\algorithmicensure}{\textbf{Output:}}
\caption{Continuous cross approximation using fiber--adaptive approximations}
\label{alg:crossapprox}
\begin{algorithmic}[1]
\Require Bivariate function $f:  [a,b] \times [c,d] \to \reals$; rank estimate $r$; initial column fiber indices $\bvec{y}=[y^{(1)}, y^{(2)}, \ldots, y^{(r)}]$; stopping tolerance $\crossdelta>0$; adaptive approximation scheme $\texttt{approx-fiber}(f_k,\approxeps)$; fiber approximation tolerance $\approxeps$
\Ensure $\bvec{x},\bvec{y}$ such that $\mat{F} \in \reals^{r \times r}$, with $\mat{F}[i,j] = f(\bvec{x}[i],\bvec{y}[j])$, has ``large'' volume
\State $\delta = \crossdelta+1$
\State $f^{(0)} = 0$
\State $k = 1$
\State Initialize columns $\mvf{C}: [a,b] \to \reals^{1 \times r}$
\State $\mvf{C}(x)[1,i] = \texttt{approx-fiber}(f(x,y^{(i)}),\approxeps)$ for $i = 1 \ldots r$ \label{alg:crossapprox:line:formC}
\While{$\delta \leq \crossdelta$}
    \State $\mvf{Q}\mat{T} = \texttt{qr}(\mvf{C})$ \Comment{QR decomposition of a matrix-valued function} \label{alg:ca:line:qr1}
    \State $\bvec{x} = \texttt{dominant}(\mvf{Q})$ 
    \State Initialize rows $\mvf{R}: [c,d] \to \reals^{r \times 1}$
    \State $\mvf{R}[j,1](y) = \texttt{approx-fiber}(f(x^{(j)}, y),\approxeps)$ for $j = 1 \ldots r$ \label{alg:crossapprox:line:formR}
    \State $\mvf{Q}\mat{T} = \texttt{qr}(\mvf{R}^T)$ \label{alg:ca:line:qr2}
    \State $\bvec{y} = \texttt{dominant}(\mvf{Q})$
    \State $\hat{\mat{Q}} = \left[\mvf{Q}[1,1](y_1)\  \mvf{Q}[1,2](y_2)\  \ldots \  \mvf{Q}[1,r](y_r)\right]$ 
    \State $\mvf{C}(x)[1,i] = \texttt{approx-fiber}(f(x,y^{(i)}),\approxeps)$ for $i = 1 \ldots r$ \label{alg:crossapprox:line:formC2}
    \State $f^{(k)}(x,y) = \mvf{C}(x)\hat{\mat{Q}}^{\dagger} \mvf{Q}^T(y)$
    \State $\delta = \Vert f^{(k)} - f^{(k-1)}\Vert  / \Vert f^{(k)}\Vert $
    \State $k = k+1$
\EndWhile
\end{algorithmic}
\end{algorithm}

Finally, we point out the distinction between this row-column alternating algorithm and \textit{adaptive} cross approximation algorithms as in~\cite{Bebendorf2000}. Adaptive cross approximation algorithms for the approximation of bivariate functions build up a set of row and column indices by sequentially choosing points which maximize the volume. In other words, they do not require one to prescribe a rank, nor do they require one to evaluate entire function fibers. These algorithms, which are equivalent to LU decomposition with complete pivoting, can use two-dimensional optimization methods to seek the optimal locations for function evaluation (which become the pivots). This methodology is attractive for bivariate problems but difficult to extend to the multivariate case, because it would require us to optimize over locations in $d$ dimensions. 

\subsection{Approximate maxvol through dominant submatrices}\label{sec:maxvolskinny}

Algorithm~\ref{alg:crossapprox} uses the \texttt{dominant} algorithm to select the next fiber locations. In particular, the goal of \texttt{dominant} is to find indices that correspond to a \textit{dominant} submatrix. For the case of matrix-valued functions, the definitions of a submatrix and a dominant submatrix are below.

\begin{definition}[Submatrix of a matrix-valued function] \label{def:sub1d}
  A submatrix of a matrix-valued function $\mvf{A}:\xspace{} \to \reals^{n \times r}$ is the matrix $\submatm{A} \in \reals^{r \times r}$ obtained by fixing a set of $r$ tuples\newline $\{(i_1,x_1), (i_2,x_2), \ldots, (i_r,x_r)\}$, where $(i_k,x_k) \in \{1,\ldots,n\} \times \xspace{}$   for $k,l = 1 \ldots r$ and $(i_k,x_k) \neq (i_l,x_l) $ for $k \neq l$, such that 
$\submatm{A}[k,j] = \mvf{A}[i_k,j](x_k).$
\end{definition}

\begin{definition}[Dominant submatrix]
A dominant submatrix of a matrix-valued function $\mvf{A}:\xspace{} \to \reals^{n \times r}$ is any submatrix $\submatm{A}$ such that for all values $(i,x,k) \in \{1,\ldots,n\} \times \xspace{} \times \{1,\ldots,r\}$ the matrix-valued function $\mvf{B} = \mvf{A}\submatm{A}^{-1}$ is bounded as $\left \vert \mvf{B}[i,k](x) \right \vert \leq 1$. \label{def:dom}
\end{definition}

Pseudocode for computing a dominant submatrix of a matrix-valued function is given in Algorithm~\ref{alg:maxvol};\footnote{We construct the algorithm for a matrix-valued function in order to generalize to multivariate cross approximation. The bivariate case is a special case where the matrix-valued function has a single row. The matrix-valued function arises in~\cite[Lines 4, 7,  8,  21 in Alg.~1]{Savostyanov2011}. In that discrete algorithm, reshapings of TT cores are formed. In our continuous case, these reshapings need not be formed. Instead they require treating each continuous core~\eqref{eq:ftcore} as a matrix whose rows are indexed by $(i_k, x_k)$ pairs.} 
it mirrors the algorithm provided in~\cite{Goreinov2010} for ``tall and skinny'' matrices. The algorithm seeks a submatrix which is {dominant}, because a dominant submatrix has a volume that is not much smaller than that of the maximum volume submatrix, \revise{and finding the actual maximum volume submatrix is computationally intractable.} For a further discussion of this approach we refer the reader to~\cite{Goreinov1997maxvol,Goreinov2001,Goreinov2010}. %

The algorithm works by swapping ``rows'' of $\mvf{A}$ (these are now specified by (index, $x$-value) combinations) until all the elements of $\mvf{B}$ are less than 1. This is exactly what the operations in Lines $\ref{alg:maxvol:opt}$ and $\ref{alg:maxvol:switch}$ of Algorithm~\ref{alg:maxvol} are doing. To find an initial set of  linearly independent ``rows'' of $\mvf{A}$, the algorithm first performs a continuous pivoted LU decomposition, denoted by \texttt{lu}, yielding pivots $\bvec{\alpha} = \{ (i_1,x_1), \ldots, (i_r,x_r) \}$.

\begin{algorithm}% [ht!]
  \caption{$\texttt{dominant}$: Dominant submatrix of a matrix-valued function}  
  \label{alg:maxvol}
\renewcommand{\algorithmicrequire}{\textbf{Input:}}
\renewcommand{\algorithmicensure}{\textbf{Output:}}
\begin{algorithmic}[1]
\Require $\mvf{A}:\xspace{} \to \reals^{n \times r}$, a matrix-valued function.
\Ensure $\bvec{\alpha} = \{ (i_1,x_1,) \ldots (i_r,x_r) \}$ such that $\submatm{A}$ is dominant
\State $\left \{ L,\mat{U},\bvec{\alpha}  \right \}= \texttt{lu}(A)$ \Comment{LU decomposition of the matrix-valued function}
\State $\delta = 2$
\While {$\delta > 1$}
    \State $\submatm{A} \leftarrow $ submatrix defined by $\bvec{\alpha}$
    \State $x^*,i^*,j^* = \underset{(x,i,j)}{\arg \max} \ \mvf{A}[i,:](x)\submatm{A}^{\dagger}[:,j]$ \label{alg:maxvol:opt}
    \State $\delta = \mvf{A}[i^*,:](x^*)\submatm{A}^{\dagger}[:,j^*]$
    \If {$\delta > 1$}
         \State $x_{j^*} = x^*$, $i_{j^*} = i^*$ \label{alg:maxvol:switch}
    \EndIf
\EndWhile
\end{algorithmic}
\end{algorithm}

There exist two differences between the continuous and discrete versions of these algorithms. The first is that we utilize a \textit{continuous} pivoted LU decomposition~\cite{Townsend2014,mythesis}. This approach allows the pivots to range freely over the input space, rather than being restricted to lie on a discretized set of locations. The second difference is the optimization problem specified in Line~\ref{alg:maxvol:opt} of Algorithm~\ref{alg:maxvol}. First, it is a \textit{continuous} optimization problem in $x$, allowing us to search over the entire space $\xspace{}$; in a tensor arising through discretization and compression in TT/QTT format, this maximization can only occur over the discretized points. If the discretization missed a local feature, then the approximation could suffer. Secondly, when $\mvf{A}:[a,b] \to \reals^{n \times r}$, the continuous optimization problem involves a one-dimensional decision variable and can exploit the structure of the scalar-valued functions that comprise $\mvf{A}$. For example, if these scalar-valued functions are represented \revise{with an expansion of} orthonormal polynomials, then the maximization reduces to an eigenvalue problem~\cite{Day2005} \revise{and therefore does not require the use of optimization approaches that might only converge to local minima}. If these functions are piecewise polynomials, then we search over each piece; if they are piecewise linear, then we search over each node. The benefit we gain is \textit{exploration} of the continuous space. 

\subsection{Rounding}

The cross approximation procedure discussed above relies on specifying the ranks $r_k$ \textit{a priori}. In this section, we discuss a procedure called \textit{rounding} that will be used to find ranks adaptively and that is useful for reducing the rank of a function formed by the multiplication and addition operations. The idea is that if a representation with certain ranks can be well approximated by a representation of smaller ranks, then we have overestimated the ranks used in the cross approximation algorithm and can be confident about its results. 

Rounding begins with a low-rank representation and aims to generate an approximation with the smallest ranks that is accurate to a specified relative error $\epsilon$. This is useful not only as a way of verifying the ranks used in cross approximation, but also because the ranks of~\eqref{eq:ft} may initially be higher than necessary, depending on how it was created. We now describe a continuous rounding procedure, similar to the procedure for discrete TT representations in~\cite[Algorithm 2]{Oseledets2011}, with the primary distinction involving the notions of left and right orthogonality. In the discrete case, these notions refer to orthogonality in certain reshapings of the three-dimensional TT cores. In our continuous case, these notions refer to orthogonality between vector-valued functions that are the rows or columns of the continuous cores.

Rounding follows directly from the definition of the ranks of the unfoldings of the function $f$. An unfolding, or $k-$\textit{separated} representation, of $f$ is a grouping of the first $k$ variables and the last $d-k$ variables to form a bivariate representation of a multivariate function. This representation is denoted by the superscript $k$
\begin{align}
f^k:\xspace{\leq k} \times \xspace{>k} \to \reals,  \quad \textrm{ such that } f^k(\{x_1,\ldots, x_k\}, \{x_{k+1},\ldots,x_{d}\}) = f(x_1,\ldots,x_d),  \label{eq:funfold}
\end{align}
where $\xspace{\leq k} = \xspace{1} \times \cdots \times \xspace{k}$ and $\xspace{>k} = \xspace{k+1} \times \cdots \times \xspace{d}$.

Suppose that we start with the first unfolding of a function
\begin{equation*}
    \funfold{1}(x_1, \foldright{1}) = \fcore{1} \left[ \fcore{2} \fcore{3} \ldots \fcore{d} \right] 
                = \fcore{1} \mvf{V}_1(x_{>1}),
\end{equation*}
where $\mvf{V}_1: \frspace{1} \to \reals^{r_1 \times 1}$. This equation expresses $\funfold{1}$ in a rank $r_1$ \revise{low rank representation. One can then use standard algorithms to compress this representation further~\cite{Bebendorf2008}; for instance we compress $\funfold{1}$ via a truncated SVD.} To this end, we first perform two QR decompositions of matrix-valued functions %
\begin{equation}
 \mvf{F}_1(x_1) = \mvf{Q}_1(x_1)\mat{R}_1 \ \  \textrm{ and }  \ \ 
 \mvf{V}_1^{\top}(\foldright{1}) = \widetilde{\mvf{Q}}_1(\foldright{1})\widetilde{\mat{R}}_1, \  \textrm{ such that } \funfold{1} = \mvf{Q}_1\mat{R}_1\widetilde{\mat{R}}_1^{\top} \widetilde{\mvf{Q}}_1^{\top}. \label{eq:roundQR}
\end{equation}
Then, we calculate the truncated SVD of $\mat{R}_1\widetilde{\mat{R}}_1^{\top} \simeq \mat{U_1D_1V_1^{\top}}$, 
where $\mat{U}_1 \in \reals^{r_1 \times \hat{r}_1}$, $\mat{V}_1^{\top} \in \reals^{\hat{r}_1 \times r_1}$, $\mat{D_1} \in \reals^{\hat{r}_1 \times \hat{r}_1}$ is a diagonal matrix, and the truncation level $\hat{r}_1 \leq r_{1}$ is chosen to obtain an approximation 
\begin{equation*}
   g_1 = \underbrace{\mvf{Q}_1\mat{U_1}}_{\hat{\fcore{1}}}  \underbrace{\left[\mat{D_1}\mat{V_1}^{\top}\widetilde{\mvf{Q}}^{\top}_1\right]}_{\hat{\fcore{2}} \hat{\fcore{3}} \ldots \hat{\fcore{d}} },
\end{equation*}
with error $\lVert f - g_1 \rVert \leq \delta$. $\hat{\fcore{1}}$ now has size $1 \times \hat{r}_1$.
Now that we have reduced the number of columns of the first core and the number of rows of the second core from $r_1$ to $\hat{r}_1$, we move on to the second unfolding of $g_1$ using the updated cores. Again, a truncated SVD is performed on this unfolding to reduce the rank from $r_2$ to $\hat{r}_2$. In this case we have
\begin{align*}
    g_1^2(\foldleft{2},\foldright{2}) &= \left[\facore{1} \facore{2}\right] \left[ \fcore{3}(x_3) \ldots \fcore{d}(x_d) \right] = \mvf{U}_2\mvf{V}_2,  \\
     \mvf{U}_2(\foldleft{2}) &= \mvf{Q}_2(\foldleft{2})\mat{R}_2, \quad \mvf{V}_2^{\top}(\foldright{2}) = \widetilde{\mvf{Q}}_2(\foldright{2})\widetilde{\mat{R}}_2
\end{align*}
where $\mvf{U}_2: \mathcal{X}_2 \to \reals^{1 \times r_2}$ is the matrix-valued QR of the left cores and $\mvf{V}_2: \mathcal{X}_2 \to \reals^{r_2 \times 1}$ and is the matrix-valued QR of the right cores. These QR decompositions are used %% of matrix-valued functions are again used
to obtain a truncated SVD:
\begin{equation}
   \mat{R}_2\widetilde{\mat{R}}_2^{\top} \simeq \mat{U_2D_2V_2^{\top}}, \qquad g_2 = \underbrace{\mvf{Q}_2\mat{U_2}}_{\facore{1}\facore{2}}  \underbrace{\left[\mat{D_2}\mat{V_2}^{\top}\widetilde{\mvf{Q}}^{\top}_2\right]}_{\facore{3} \ldots \facore{d}},
\end{equation}
where the first equation is the truncated SVD and the second is the new approximation
with $\lVert g_2 - g_1 \rVert \leq \delta$. 
After the truncation associated with the first and second cores, we obtain a total error of
$\lVert g_2 - f \rVert = \lVert g_2 - g_1 + g_1 - f \rVert \leq \delta + \delta  = 2\delta.$
We repeat this procedure $(d-1)$ times to obtain a final approximation $\hat{f} = g_{d-1}$ such that 
$\lVert \hat{f} - f \rVert \leq (d-1) \delta$ and
a relative error $\epsilon$ by setting $\delta = \frac{\epsilon}{d-1}\lVert f \rVert$.

The main computational burden of the algorithm described above is calculating the QR decompositions of the matrix-valued functions $\mvf{U}_k(\foldleft{k})$ and $\mvf{V}_k(\foldright{k})$, since these functions have potentially high-dimensional inputs. A computationally feasible rounding procedure requires that these QR decompositions be tractable. In the discrete setting, one can obtain an algorithm that only requires the QR decompositions of reshapings of single cores. In the continuous setting, we can similarly construct an algorithm that requires only the QR decomposition of each univariate core. This algorithm starts by assuming that all the cores $\fcore{2}, \ldots, \fcore{d}$ have \textit{orthonormal rows}. Then, we can show that $\mvf{V}_1$ also has orthonormal rows and that we are not required to take its QR decomposition in~\eqref{eq:roundQR}. Thus, in the first step of the rounding procedure we only need to compute the QR decomposition of $\mvf{F}_1$ \eqref{eq:roundQR}. The notion of \textit{orthonormal rows}, 
that $\left  \langle \fcore{k}[i,:], \fcore{k}[j,:] \right  \rangle = \delta_{i,j}$  for $i,j = 1 \ldots r_{k-1}$, 
is analogous to left orthogonality~\cite{Oseledets2011}; while \textit{orthonormal columns},
that $  \left  \langle \fcore{k}[:,i], \fcore{k}[:,j] \right \rangle =  \delta_{i,j}$
for $i,j = 1 \ldots r_{k}.$,
are analogous to right orthogonality.
Now that the continuous notions of orthonormality is clear, %% Using these notions of orthonormal rows and columns, one can apply
Algorithm 2 from~\cite[Algorithm 2]{Oseledets2011} can be applied by replacing 
the discrete QR decompositions with continuous QR decompositions of matrix-valued functions~\cite{mythesis}.

Overall, constructing \revise{a} low-rank approximation of a black-box function via cross approximation and rank adaptation can be implemented by successively increasing or ``kicking'' the TT ranks until rounding leads to a reduction of all ranks. This approach ensures that the ranks are \textit{overestimated} for the cross approximation procedure. Pseudocode for rank-adaptation is given in Algorithm~\ref{alg:rankadapt}.
\begin{algorithm}% [ht!]
    \caption{$\texttt{ft-rankadapt}$: FT approximation with rank adaptation}
    \label{alg:rankadapt}
\renewcommand{\algorithmicrequire}{\textbf{Input:}}
\renewcommand{\algorithmicensure}{\textbf{Output:}}
\begin{algorithmic}[1]
\Require A multivariate function $f:\mathcal{X} \to \reals$, with $\mathcal{X} \subset \reals^d$; cross approximation tolerance $\delta_{\text{cross}}$; size of rank increase $\texttt{kickrank}$; rounding accuracy $\epsilon_{\text{round}}$; initial ranks $\bvec{r} = (1, r_1, \ldots, r_{d-1}, 1)$
\Ensure Low-rank approximation $\hat{f}$ such that rank increase followed by rounding does not change ranks
\State $\hat{f} = \texttt{cross-approx}(f,\bvec{r},\delta_{\text{cross}})$ 
\State $\hat{f}_r = \texttt{ft-round}(\hat{f},\epsilon_{\text{round}})$
\State $\hat{\bvec{r}} = \textrm{ranks}(\hat{f}_r)$
\While {$\exists \, i \textrm{\ s.t.\ } \hat{r}_i = r_i$} %\hat{\bvec{r}} \neq \bvec{r}$}
    \For{$k=1$ {\bf to} $d-1$}
        \State $\bvec{r}_k = \hat{\bvec{r}}_k + \texttt{kickrank}$
        \State $\hat{f} = \texttt{cross-approx}(f,\bvec{r},\delta_{\text{cross}})$ 
    \EndFor
    \State $\hat{f}_r = \texttt{ft-round}(\hat{f},\epsilon_{\text{round}})$
    \State $\hat{\bvec{r}} = \textrm{ranks}(\hat{f}_r)$
\EndWhile
\State $\hat{f} = \hat{f}_r$
\end{algorithmic}
\end{algorithm}
\section{Numerical experiments}\label{sec:numerics}

This section describe several numerical experiments that explore the capabilities and performance of our computational approach. In Section~\ref{sec:gaussbump} we %%
demonstrate how our adaptive approach can outperform, by several orders of magnitude, the QTT and the tensor-of-coefficients approaches for approximating functions with local features. 
 In Section~\ref{sec:fiberadapt}, we show two examples of how our fiber-adaptive continuous algorithms explore a function, i.e., \emph{where} they place function evaluations in the domain. These examples offer an intuitive explanation of the observed performance gains. In Section~\ref{sec:adaptintdiff} we demonstrate that the FT can be used for adaptive integration \textit{and} differentiation, even with functions that contain a discontinuity. Furthermore, these operations use the functional representation to avoid constructing discretizations that are specific to either differentiation or integration. %%
In Section~\ref{sec:appapprox}, using a prototypical application in uncertainty quantification, we study how the fiber adaptation tolerance interacts with rounding tolerance.
All of the experiments described below are performed with the Compressed Continuous Computation $(C^3)$ library available at \url{http://github.com/goroda/Compressed-Continuous-Computation}.

We briefly mention how to perform integration and differentation in low-rank FT format. Integration involves integrating all of the univariate functions in each core and then performing matrix-vector multiplication $d-1$ times:
\begin{align*}
    \int f(x) dx &= \int \fcore{1} \fcore{2} \cdots \fcore{d} dx_1 \ldots dx_d \nonumber  \\
                &= \left(\int \fcore{1}dx_1 \right) \left( \int \fcore{2}dx_2\right) \cdots \left(\int \fcore{d} dx_d \right) 
                = \mat{\Gamma_1} \mat{\Gamma_2} \cdots \mat{\Gamma_d},  \label{eq:integrate}
\end{align*}
where $\mat{\Gamma_k} = \int \fcore{k}dx_k$ contains entries $\mat{\Gamma_k}[i,j] = \int \ffiber{k}{ij}(x_k) dx_k.$ Furthermore, since each univariate function is represented in a continuous format,
this integral is uniquely defined and computationally inexpensive. Differentiation requires differentiating the scalar-valued functions that make up the corresponding core. For example, consider the partial derivative of a $d$--variate function:
\begin{equation*}
  \frac{\partial f}{\partial x_k} = \mvf{F}_1\cdots\mvf{F}_{k-1} \left[ 
    \begin{array}{ccc}   
       \frac{d \ffiber{k}{11}}{dx_k} & \cdots & \frac{d \ffiber{k}{1r_{k}}}{dx_k} \\
        \vdots & & \vdots \\
       \frac{d \ffiber{k}{r_{k-1}1}}{dx_k} & \cdots & \frac{d\ffiber{k}{r_{k-1}r_{k}}}{dx_k}
    \end{array}
   \right] \mvf{F}_{k+1} \cdots \mvf{F}_d.
\end{equation*}
Again, because each univariate function is represented in a continuous format (e.g., in its own orthogonal polynomial basis), this operation is uniquely defined and computationally inexpensive. %%

\subsection{Approximation of local features in analytic functions}\label{sec:gaussbump}

In this section we compare local fiber-based adaptivity to discrete tensor-based approaches. In particular, we compare our fiber-adaptive scheme with the spectral tensor-train (STT) algorithms found in~\cite{Bigoni2016}. The STT approach works by first discretizing a function on a tensor-product grid defined by Gaussian quadrature, performing a low-rank decomposition of this discretized function using the TT or QTT, and projecting this decomposition onto a polynomial basis in order to obtain \textit{spectrally convergent} approximations for analytic functions. We will demonstrate that even for \textit{analytic} functions, a fiber adaptive approach can be more efficient compared with the ``discretized function'' \textit{and} ``tensor-of-coefficient'' approaches.
Our demonstration repeats an experiment from~\cite{Bigoni2016}\footnote{We use the TensorToolbox package (\url{https://pypi.python.org/pypi/TensorToolbox/}) for the STT results.} on an \textit{analytic} function that can be well represented using a global polynomial basis.
The goal is to integrate a Gaussian bump\
\begin{equation}\label{eq:gaussbump}
  f(x) = A\exp\left( -\frac{(x_1-c)^2 + (x_2-c)^2 + (x_3-c)^2}{2l^2}\right),\ \  x \in [0,1]^3, \quad  \textrm{ with } c=0.2 \textrm{ and } l=0.05,
\end{equation}
with $A=1$. As part of its construction, the STT uses a QTT to perform \textit{rank-revealing} interpolation of the function on the tensor-product grid of Gaussian quadrature nodes using the \texttt{TT-DMRG-cross} algorithm with $\epsilon=10^{-10}$. \cite{Bigoni2016} notes that QTT is more efficient than TT because it is able to localize function evaluations around the bump.

For this comparison, our \textit{continuous} algorithms employ a fiber-adaptive scheme based on piecewise polynomials with  $\crossdelta=\roundeps=10^{-10}.$ In particular, for each fiber encountered through cross approximation in Algorithm~\ref{alg:rankadapt}, the following steps are taken:
\begin{enumerate}
  \item Approximate the fiber with a degree--7 Legendre polynomial basis.
  \item If the square of \revise{the} coefficient of the highest-degree polynomial, normalized by the squared norm of the function, is greater than a tolerance $\approxeps$, then split the domain into three regions.
  \item If a split occured, return to step 1 for each new subdomain; otherwise stop.
\end{enumerate}

\begin{figure}
  \centering
  \includegraphics[scale=0.65]{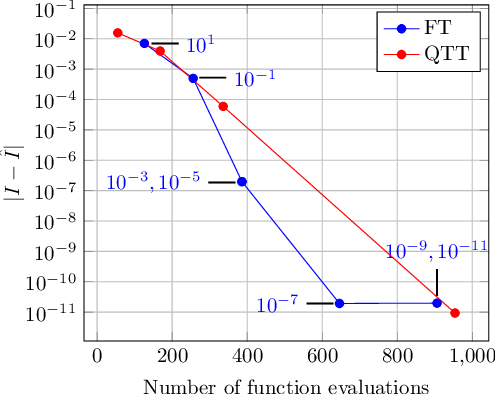}
  \caption{Convergence comparison between the spectral QTT~\cite{Bigoni2016} and the fiber-adaptive FT for a three-dimensional Gaussian bump. In the FT, each fiber is adapted using a breadth-first splitting of piecewise polynomials into three equal regions until the highest-degree coefficient in each region is below the tolerance indicated next to the data points. Rank adaptation uses \texttt{ft-rankadapt}.}
  \label{fig:gberr}
\end{figure}

The resulting errors in the integral for a given number of function evaluations are shown in Figure~\ref{fig:gberr}. Since the STT is not an adaptive scheme, increasing numbers of function evaluations are chosen by increasing the number of quadrature nodes in each dimension. The polynomial orders for the STT are increased such that the number of quadrature nodes grows in powers of two, from 4 to 32, to be conducive to the QTT. In contrast, the continuous cross-approximation algorithms we describe here are adaptive, and therefore increasing numbers of evaluations are obtained by tightening the fiber approximation tolerance $\approxeps$; tolerance values corresponding to each data point are indicated on the plot.
Both of these approximations are able to eventually achieve $\mathcal{O}(10^{-11})$ error, but the QTT-based scheme reaches this error with a larger number of evaluations. This suggests that the base-2 reshapings of the QTT do not exploit local structure as effectively as our chosen piecewise approximation. Second, note that an adaptive tolerance allows for more granular control of the number of evaluations rather than a global polynomial order. Finally, we are able to achieve a maximum of \textit{3 orders of magnitude} reduction in error for a specified number of function evaluations.

\begin{figure}
  \centering
  \begin{subfigure}[b]{0.3\textwidth}
    \includegraphics[width=\textwidth]{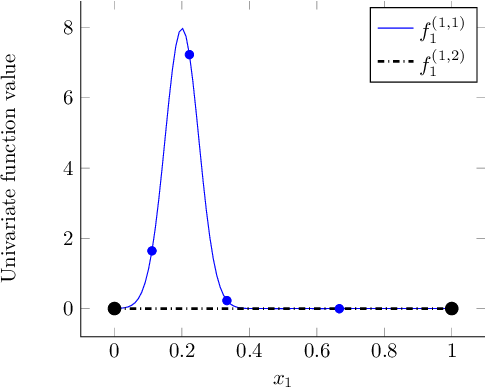}
    \caption{First core}
    \label{fig:gpbw:a}
  \end{subfigure}
  \begin{subfigure}[b]{0.3\textwidth}
    \includegraphics[width=.9\textwidth]{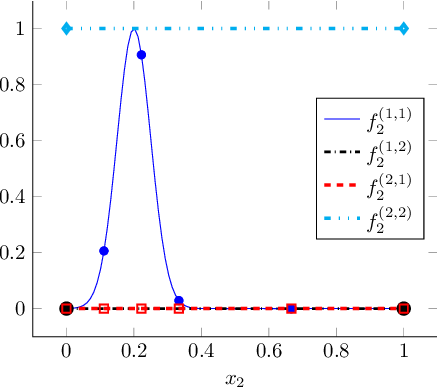}
    \caption{Second core}
    \label{fig:gpbw:b}
  \end{subfigure}
  \begin{subfigure}[b]{0.3\textwidth}
    \includegraphics[width=.9\textwidth]{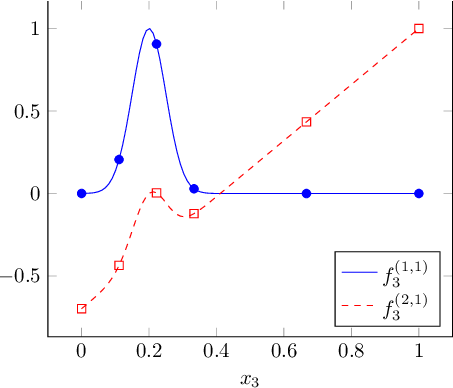}
    \caption{Third core}
    \label{fig:gpbw:c}
  \end{subfigure}
  \caption{Univariate functions for a fixed TT-rank=2 approximation. Marks indicate the location of knots found through adaptive piecewise polynomial procedures. Note that in the first and second core, there exist univariate functions with differing \textit{numbers} of knots. In other words, the basis for each function was chosen \textit{online}, \textit{adaptively}, and \textit{independently} for each fiber.}
  \label{fig:gbpw}
  \vspace{-10pt}
\end{figure}

Next we explore \textit{how} the automated adaptation algorithm chooses to approximate fibers. In particular, we perform cross-approximation with a fixed rank-2 approximation. Since~\eqref{eq:gaussbump} is actually rank-1, an efficient rank-2 approximation should introduce \textit{zero} functions that require little to no storage, as described in Section~\ref{sec:datastruct}. In Figure~\ref{fig:gbpw}, we see that this is indeed the case. To demonstrate this fact consider that~\eqref{eq:gaussbump} can be written as $f(x) = e_1(x_1)e_2(x_2)e_3(x_2)$ with $e_s$ proportional to a  univariate Gaussian bump centered at $x_j=0.2$. A TT-rank 2 function could then take the form\footnote{One can specify a different combination of zero, nonzero, and arbitrary functions by rearranging the order of the first core, for example. We specify this combination for illustrative purposes because it is the one chosen by the cross-approximation algorithm.}
\begin{equation}
f(x_1,x_2,x_3) = \left[ e_1(x_1) \ 0 \right] \left[ 
  \begin{array}{cc}
    e_2(x_2) & 0 \\
    \ffiber{2}{21}(x_2) & \ffiber{2}{22}(x_2)
  \end{array} \right]
\left[ \begin{array}{c} e_3(x_3) \\ \ffiber{3}{21}(x_3) \end{array} \right],
\end{equation}
where $\ffiber{2}{12}$, $\ffiber{2}{22}$ and $\ffiber{3}{21}$ can be arbitrary functions since they are multiplied by the zero function in the first core.

Figure~\ref{fig:gbpw} indicates that this is exactly the structure found using continuous cross approximation. The first core essentially specifies $\ffiber{1}{11}(x_1) = e_1(x_1)$ and $\ffiber{2}{12}=0$. Also notice that $\ffiber{1}{11}$ was refined twice, while $\ffiber{2}{12}$ was never refined. This adaptation happened in an \textit{online} fashion; it would have been difficult \textit{a priori} to create a basis for the first dimension that included this precise refinement. The second core also exhibits the structure we expect: $\ffiber{2}{11}(x_2) = e_2(x_2)$, $\ffiber{2}{12}(x_2) = 0$, $\ffiber{2}{21}(x_2)=0$, and $\ffiber{2}{22}(x_2) = 1$. Note that the latter two functions %% 
could have been chosen \textit{arbitrarily}. In this case, they were chosen to be constants. The third core chooses $\ffiber{3}{11} = e_3(x_3)$ and an arbitrary function for $\ffiber{3}{21}$.

As described above, in this example we chose a particular uniform refinement scheme for piecewise polynomials. Other choices of number of regions or polynomial degree will naturally lead to different performance characteristics. Our goal was simply to show that the flexibility of \textit{allowing} fiber adaptation \textit{can} greatly improve performance. More complex non-uniform and function-adaptive refinement schemes~\cite{Cohen2009} for piecewise polynomials can certainly be used. Our approach allows for \textit{any} univariate function approximation scheme. In practice, this scheme can and should be tailored to any problem-specific knowledge.

\begin{figure}
  \centering
  \begin{subfigure}[b]{0.3\textwidth}
    \includegraphics[width=\textwidth]{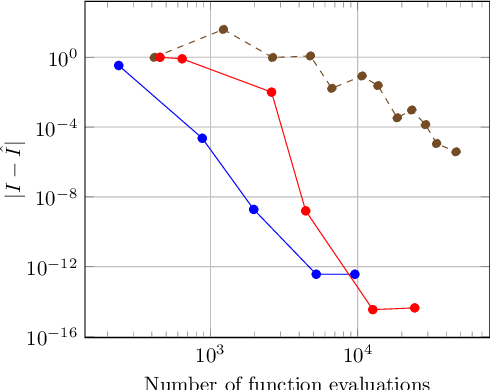}
    \caption{$d=8$}
    \label{fig:gd-8}
  \end{subfigure} 
  \begin{subfigure}[b]{0.3\textwidth}
    \includegraphics[width=\textwidth]{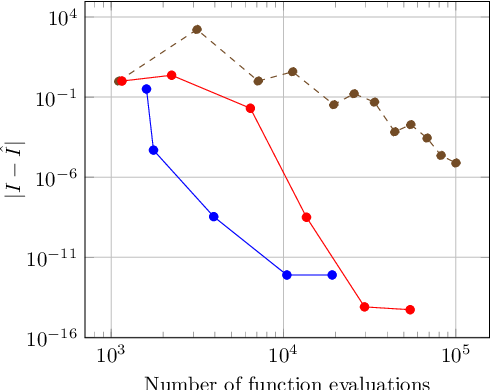}
    \caption{$d=16$}
    \label{fig:gd-16}
  \end{subfigure}
  \begin{subfigure}[b]{0.3\textwidth}
    \includegraphics[width=\textwidth]{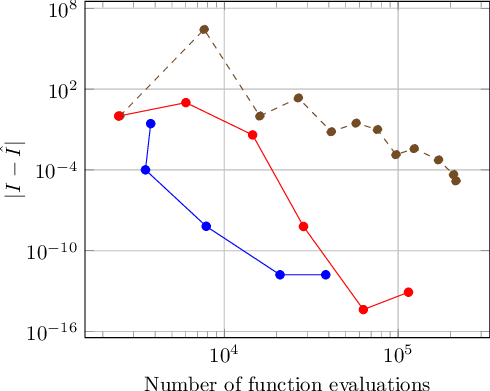}
    \caption{$d=32$}
    \label{fig:gd-32}
  \end{subfigure}
  \caption{Comparison of convergence for FT (blue), QTT (red), and TT (brown) for the probability density functions of $d$-dimensional normal random variables. FT uses fiber-adaptive cross approximation, while the QTT and TT represent the coefficients of global polynomials and are computed using DMRG-cross. }
  \label{fig:gaussdim}
\end{figure}

Finally, we compare the fiber-adaptive FT %% 
with the QTT- and TT-based spectral decompositions for approximating the probability density function (PDF) of normal random variables of increasing dimension. Our tests are performed on $d$-dimensional versions of \eqref{eq:gaussbump}. %% 
We refine the spectral approximations by increasing polynomial degree, and, correspondingly, we tighten the univariate approximation tolerance for the FT. The results for the error in the integral are shown in Figure~\ref{fig:gaussdim}. Again, we observe excellent performance of the fiber-adaptive FT scheme. 

\subsection{Fiber adaptivity examples}\label{sec:fiberadapt}

We now explicitly demonstrate that the locations of function evaluations resulting from fiber adaptation do not lie on a tensor product grid; instead they are adapted to the local features.
 These examples illustrate one of the primary ways that our approach differs from the discrete tensor-train and from other functional tensor-train approaches. Consider the approximation of two canonical rank-one functions:
\begin{align}
f_{\text{sin}}(x_1,x_2) &= \sin\left(10x_1+ 1/4\right)\left(x_2+1\right), \qquad (x_1,x_2) \in [-1,1]^2 \label{eq:sin} \\
f_{\text{genz2d}}(x_1,x_2) &= \left\{
  \begin{array}{cc}
    0 & \textrm{ if } x_1 > 0.5 \textrm{ or } x_2 > 0.5 \\
    \exp(5x_1 + 5x_2) & \textrm{ otherwise }
  \end{array} \right., \qquad (x_1,x_2) \in [0,1]^2 \label{eq:2dgenz}
\end{align}
where~\eqref{eq:2dgenz} is a bivariate Genz function of the `discontinuous' family~\cite{Genz1984}.

For~\eqref{eq:sin} we use a global expansion of Legendre polynomials 
to represent each fiber. %
The coefficients are determined via projection using Clenshaw-Curtis quadrature \revise{to obtain nested rules; other quadrature rules, for instance Gaussian quadrature, can also be readily used}. The degree $n-1$ of each fiber approximation is progressively increased until four successive coefficients %
 have magnitudes less than $\approxeps = 10^{-10}$. 

\begin{figure}
\begin{center}
  \begin{subfigure}[b]{0.43\textwidth}
    \includegraphics[width=0.85\textwidth]{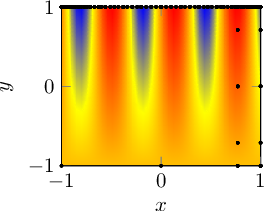}
    \caption{Sine function~\eqref{eq:sin} and evaluation points.}
    \label{fig:2d:a}
  \end{subfigure}
  \begin{subfigure}[b]{0.43\textwidth}
    \includegraphics[width=0.85\textwidth]{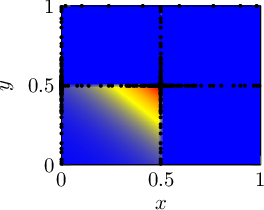}
    \caption{Genz function~\eqref{eq:2dgenz} and evaluation points.}
    \label{fig:2d:b}
  \end{subfigure}
\caption{Contour plots and evaluations of $f_{\text{sin}}$ and $f_{\text{genz2d}}$.}
\label{fig:2dexamples}
\end{center}
\end{figure}

Figure~\ref{fig:2dexamples}(a) shows the function and the parameter values where it is evaluated by the cross approximation algorithm. Even though~\eqref{eq:sin} has rank one, we deliberately construct a rank-two cross approximation to illuminate the pattern of function evaluations chosen by our method. We observe that fibers at different positions are approximated using differing numbers of function evaluations. The number of evaluations required in the oscillatory region (near $x_2=1$) is much greater than in the constant portion ($x_2=-1)$. Such an adaptation of the grid would be difficult to achieve with discrete TT, and highlights the flexibility of the continuous approximation framework.

For~\eqref{eq:2dgenz} we can no longer use a global polynomial expansion due to the discontinuity, and therefore we employ piecewise polynomial approximations. We perform adaptation in a similar manner to  Section~\ref{sec:gaussbump}. We use a breadth-first refinement where we split polynomials into three regions, approximate each region with a degree--6 polynomial, and further split the region if the leading coefficient of the expansion %% 
contributes more than $\approxeps=10^{-10}$ to the squared norm of the univariate function. % 
We then approximate each smooth interval of the fiber using suitably scaled Legendre polynomials. In this example, we seek a rank-one approximation of the function. Function contours and evaluation locations are shown in Figure~\ref{fig:2dexamples}(b). This approximation, like that of~\eqref{eq:sin}, achieves machine accuracy. We see that the algorithm clusters evaluations points around the discontinuity, as desired, and that the evaluations again do not lie on a tensor-product grid.

Having constructed low-rank representation of these functions, we can perform computations directly in compressed format. For example, we can now integrate the discontinuous function $f_{\text{genz2d}}$. Such an integration could not be performed using array-based tensor-train algorithms, unless the entire domain were manually partitioned; otherwise, one would need specialized integration rules to deal with the discontinuity. By representing everything in functional form, we are able to perform integration and approximation automatically.
In Section~\ref{sec:genzint}, we will evaluate integration performance for discontinuous Genz functions on higher-dimensional input spaces.

\subsection{Adaptive integration and differentiation}\label{sec:adaptintdiff}

In this section, we fix a target fiber approximation tolerance and fix the ranks to
determine the resulting errors of approximate integrals and derivatives. We show that the error is well behaved and sufficient for most applications. Our second goal is to show that our representation of the FT is flexible enough to enable a variety of different computations. Specifically, we demonstrate that we simultaneously obtain adaptive integration and differentiation procedures by using an adaptive approximation approach.
No specialized integration or finite difference algorithms or grids are required.

\subsubsection{Integration of rank--two sine function}\label{sec:sinint}
As in~\cite{Oseledets2010} we consider the rank-two function
\begin{equation}
	f(x) = \sin(x_1 + x_2 + \ldots + x_d),
	\label{eq:sinsum}
\end{equation} 
whose integral on the unit hypercube is known analytically for any $d$:
$\int_{[0,1]^d} f(x) dx = \textrm{Im}\left[ \left(  \frac{e^i - 1}{i} \right)^d \right].$

We study the performance of FT-based integration as a function of dimension and the fiber adaptation tolerance. Specifically, we approximate each univariate fiber with a Legendre polynomial expansion using pseudospectral projection with Gaussian quadrature. Each fiber is initially represented by an expansion of degree $k=5$, and the degree is increased from $k$ to $k+7$ with each adaptation step. We stop adaptation after the last \textit{two} coefficients contribute less than $\approxeps$ to the squared norm of the function. 
\begin{figure}% 
\begin{center}
\includegraphics[scale=0.8]{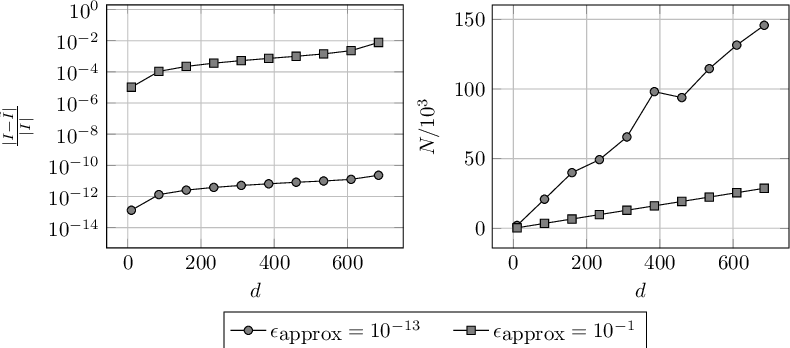}
\caption{Relative errors (left panel) and growth in the number of evaluations (right panel) involved in the integration of~\eqref{eq:sinsum}, as a function of dimension $d$ and fiber adaptation parameter $\approxeps$.}
\label{fig:sinsum}
\end{center}
\end{figure}
Figure~\ref{fig:sinsum} shows the relative error in the integral as a function of $\approxeps$ and the dimension $d$ of the problem. The right panel also shows the number of function evaluations used to construct the FT approximation. As desired, the number of evaluations grows \textit{linearly} with the input dimension, and %% 
the integration error is well behaved across dimensions, even for $d > 600$. (Note that the error axis is on a log scale and hence the variation of the errors for the tightest tolerance would be virtually invisible on the $\approxeps = 10^{-1}$ curve.) %% 

\subsubsection{Integration and differentiation of rank--one discontinuous Genz function}\label{sec:genzint}

We now demonstrate integration and differentiation of the discontinuous Genz function for varying dimensions. Specifically, we consider $f: [0,1]^d \to \reals$ defined as 
\begin{equation}
f(x_1,x_2,\ldots,x_d) = \left\{
\begin{array}{cc}
0 & \textrm{ if } x_i > \frac{1}{2} \textrm{ for any } i = 1 \ldots d \\
\exp\left(\sum_{i=1}^{d} 5x_i\right)  & \textrm{ otherwise }
\end{array} . \right.
\label{eq:discgenz}
\end{equation}
The exact integral of~\eqref{eq:discgenz} is
$
I[f] = \left( \frac{\exp\left(\frac{5}{2}\right) - 1}{5}\right)^d
$;
this problem is quite challenging because the integral grows exponentially with dimension. For example, for $d=10$ the value of the integral is $I\approx 3.131 \times 10^{3}$, while for $d=100$ the value is $I \approx 9.05455 \times 10^{34}$. In addition to computing the integral, we wish to evaluate numerically the following derivative:
$
D[f](x) = \nabla \cdot f(x_1,\ldots,x_d) = 5d\exp(\sum_{i=1}^d5x_i)
$
for $x_i = 0.2$ for all $i$.

Fiber adaptation is performed using the same scheme as for the two-dimensional case of Section~\ref{sec:fiberadapt}. For this function, individual fiber adaptation is essential for approximating in the presence of discontinuities. \revise{Note that other discontinuity detection schemes~\cite{Archibald2005,Archibald2009,Gorodetsky2014} could also be readily used within our framework by embedding them within univariate function approximation routines.}

\begin{figure}% [ht!]
\begin{center}
\includegraphics[width=0.8\textwidth]{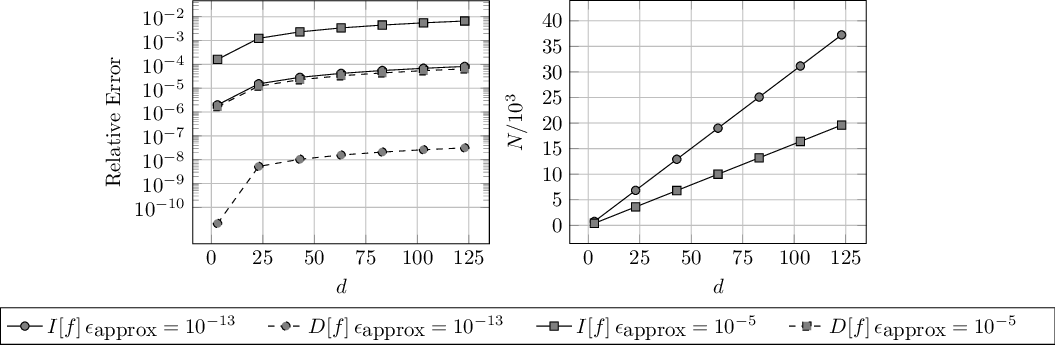}
\caption{Relative errors (left panel) and \# evaluations (right panel) when integrating or differentiating~\eqref{eq:discgenz}, versus dimension $d$, for different values of the fiber approximation tolerance $\approxeps$.}
\label{fig:genznd}
\end{center}
\end{figure}

The resulting relative errors for $I$ and $D$, and the required number of function evaluations, are shown in Figure~\ref{fig:genznd} as functions of the dimension $d$.

The results indicate stable approximations for all values of $d$ for fixed tolerance levels, and that the number of function evaluations scales linearly with dimension. Furthermore, we are able to approximate extremely large values of the integral, suggesting a general robustness of the algorithm. 

Simultaneous integration and differentiation as reported here would be extremely difficult to perform using either the discrete tensor-train or the spectral tensor-train~\cite{Bigoni2016} techniques, because discontinuities pose problems for most integration rules. For example, we can compare our adaptive integration scheme to a scheme based on finite differences using a QTT decomposition of a function discretized on a regular grid. Practically, we can consider discretizations as fine as $2^{50}$ in each input coordinate, because the complexity of the QTT grows like the log of the discretization size. However, such discretizations can be difficult to employ in an adapt-to-tolerance setting, due to the difficulty of controlling roundoff errors. %% 

\begin{figure}% [ht!]
\begin{center}
\includegraphics[width=0.35\textwidth,angle=90]{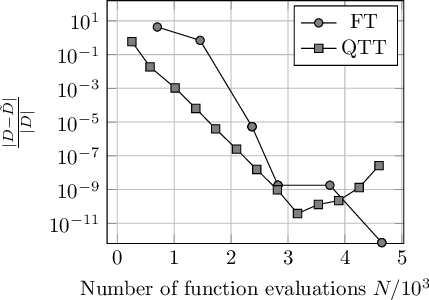}
\caption{Relative error of differentiation with the continuous tensor-train using adaptive piecewise polynomials and the QTT using finite differences for a ten dimensional version of~\eqref{eq:discgenz}. FT refinement is obtained by diminishing $\approxeps$. QTT refinement is performed by quadrupling the size of the grid along each dimension. Refining functions proves more numericaly stable than refining grids.}
\label{fig:fd-comparison}
\end{center}
\end{figure}

Figure~\ref{fig:fd-comparison} shows, for $d=10$, the error of the derivative with increasing numbers of function evaluations. The number of evaluations required by the QTT is governed by the size of discretization, while the number of evaluations for the FT is governed by the fiber adaptation tolerance. Errors with the QTT approach decrease rapidly, and then increase just as sharply. Our FT approach is more robust in the sense that it avoids such numerical instabilities (which are due to roundoff error). In effect, increasing the number of function evaluations by tightening a tolerance is safer than simply refining a tensor-product discretization. Practically, it is difficult to decide how finely to discretize a function and then to assess the accuracy of a fixed discretization scheme; our adaptive approach instead learns more about the structure of the problem in order to avoid such choices.

\revise{
\subsection{Approximation of a function with a singularity}\label{sec:singularity}

In this section we explore the properties of our algorithm on a function that does not readily admit low-rank structure and that has a singularity (which we soften by adding a small value to the denominator), namely
\begin{equation}\label{eq:singularity}
  f(x_1,\ldots,x_d) = \frac{1}{\sqrt{\sum_{i=1}^d x_i^2 + 10^{-12}}},
\end{equation}
for $x \in [-1,1]^d.$ The quadratic in the denominator of this function is a sum of univariate functions and is therefore rank 2 regardless of dimension. However, taking its square root and then inverse causes the rank to grow with dimension. We seek to study how this rank grows with dimension and how the rank adaptation scheme performs with varying tolerances.

For this problem we again use piecewise-polynomials of maximum degree three, splitting the domain of a univariate fiber into four regions for adaptation. We initialize the rank to four, use a maximum of five cross-approximation iterations, and set $\texttt{kickrank}=2.$ We report the root relative mean squared error of the approximation on $10,000$ uniformly sampled inputs, and we provide the average FT rank found by the adaptive rouding scheme. These results are shown in Figure~\ref{fig:sing}. The left panel shows the approximation error as a function of dimension for various rounding tolerances. First, the approximation error decreases with dimension but is still generally larger than the errors reported elsewhere in this paper. This fact reflects the challenging nature of this target function. The reduction in error with dimension can be attributed to the fact that the relative volume of the space that is affected by the discontinuity actually \textit{decreases} with dimension. Furthermore, we see that over the given range of rounding tolerances, there is no significant difference in the approximation accuracy. This can be attributed to the fact that univariate approximation errors are the dominant source of inaccuracy, due to the difficulty of approximating the function near the origin. In Section~\ref{sec:appapprox} we will provide an example where the opposite occurs.

\begin{figure}
  \centering
  \begin{subfigure}[b]{0.48\textwidth}
    \centering
    \includegraphics[width=\textwidth]{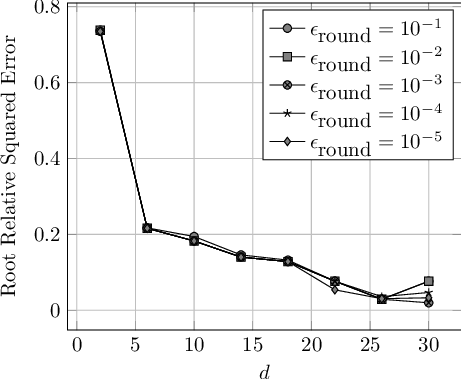}
    \caption{Approximation Errors}
    \label{fig:sing:a}
  \end{subfigure}
  \begin{subfigure}[b]{0.48\textwidth}
    \centering
    \includegraphics[width=\textwidth]{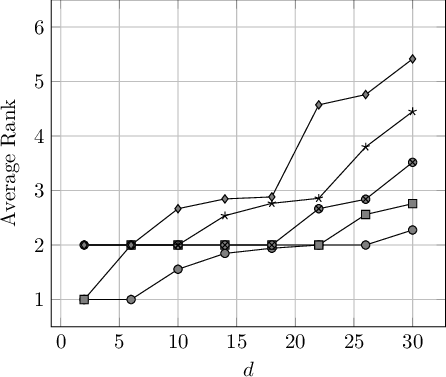}
    \caption{Average ranks}
    \label{fig:sing:b}
  \end{subfigure}
  \caption{Approximation errors and average rank for various $\epsilon_{round}$ on the function with a singularity~\eqref{eq:singularity}.}
  \label{fig:sing}
\end{figure}

The right panel of Figure~\ref{fig:sing} shows the average ranks determined by the rank adaptation scheme. As expected the rank grows with dimension. Since the errors with dimension beyond five are relatively steady, this plot indicates the correct behavior that rank needs to grow to maintain a given level accuracy. Furthermore, the tighter the rounding tolerance, the larger the average rank. To demonstrate the process of adaptivity, in Table~\ref{tab:adaptivity} we show the ranks following the rounding portion of each stage of Algorithm~\ref{alg:rankadapt}.

\begin{table}
  \centering
  \begin{tabular}{|c|c|}
    \hline
    Iteration & Ranks \\
    \hline
    \hline
    1 & 1 2 4 4 4 4 4 4 4 4 4 4 4 4 4 4 4 4 4 4 4 4 4 4 4 4 4 3 3 3 1 \\
    2 & 1 2 4 5 6 6 5 6 6 5 5 5 5 5 5 5 6 6 5 5 5 5 5 5 5 5 5 4 4 3 1 \\
    3 & 1 2 4 5 6 6 5 5 6 5 5 5 6 5 5 5 5 6 5 6 5 5 5 5 5 5 5 5 5 4 1 \\
    4 & 1 2 4 6 6 6 6 6 7 6 6 5 6 5 5 5 5 5 5 5 6 6 6 6 6 6 6 5 5 5 1 \\
    5 & 1 2 4 5 6 7 7 7 7 7 6 6 5 5 5 5 5 6 5 5 6 6 6 5 5 5 5 5 5 5 1 \\
    6 & 1 2 4 5 6 7 7 7 7 7 7 7 6 5 5 5 5 5 5 5 5 5 5 5 5 5 5 5 5 5 1 \\
    \hline
  \end{tabular}
  \caption{Function ranks after rounding for given iteration of Algorithm~\ref{alg:rankadapt} for Equation~\eqref{eq:singularity} with $\epsilon_{round} = 10^{-5}$ and $d=30$.}
  \label{tab:adaptivity}
\end{table}

}

\subsection{Approximation of an elliptic PDE}\label{sec:appapprox}
We next explore the effects of various parameters of the FT approximation algorithm on a model of subsurface flow frequently encountered in UQ applications.
 Consider the following one-dimensional elliptic PDE:
$\frac{\partial}{\partial s}\left( k(s,\omega) \frac{\partial u}{\partial s}\right) = s^2,$
for $s \in [0,1]$ with $u(s)|_{s=0} = 0$ and $\frac{\partial u}{\partial s}|_{s=1} = 0$. We consider the effects of an uncertain permeability field $k$ on a functional of the output pressure $u$. The permeability is modeled as a random process $k(s,\omega)$, endowed with a log-normal distribution:
$
\log\left[k(s,\omega) - a\right] \sim \mathcal{N}(0,c(s,s^{\prime})).
$
Here the covariance kernel is chosen to be  $c(s,s^{\prime}) = \sigma^2 \exp \left ( - \frac{|s-s^{\prime}|}{l} \right)$. To obtain a finite dimensional representation of $k(s,\omega)$, we use the Karhunen-Lo\`{e}ve expansion of $\log \left ( k(s, \omega) - a \right )$ to express the random field using the eigenfunctions of an integral operator with kernel $c$. In particular, given eigenfunctions and eigenvalues obtained from $\int_{[0,1]} c(s,s^{\prime}) \phi_i(s^{\prime}) ds^{\prime} = \lambda_i \phi_i(s)$, the log-normal process is represented as $k(s,\omega) = a + \exp \left (\sum_{i=1}^{\infty}\sqrt{\lambda_i} \phi_i(s) \xi_i(\omega) \right )$, where $\xi_i(\omega)$ are independent standard Gaussian random variables. We truncate this expansion after 24 modes.

The approximation objective of this problem is the solution of the PDE at a spatial location $s$ for a realization of the permeability field parameterized by $(\xi_i)_{i=1}^{24}$. For simplicity, we will fix $s=0.7$ to obtain a quantity of interest $Q$ that is only a function of the $\xi_i$, i.e., $Q := u(0.7,\xi_{1},\ldots,\xi_{24})$.   
We measure the relative $L^2$ error of our approximation with $n=10^4$ Monte Carlo samples as
\begin{equation*}
\text{error} = \sqrt{ \frac{\sum_{i=1}^{n} (Q(\xi) - \widehat{f}(\xi))^2}{\sum_{i=1}^{n} Q(\xi)^2}}.
\end{equation*}
We will construct approximations of $Q$ using  parameter settings for $(a,\sigma^2,l)$ that correspond to three different levels of difficulty for this problem. In particular, we will tackle an ``easy'' problem (P1) where  $(a,\sigma^2,l) = (0.0,0.1,0.125)$, a ``moderately'' difficult problem  (P2) where  $(a,\sigma^2,l) = (0.5,1.0,0.045)$, and a ``harder'' problem (P3) where $(a,\sigma^2,l) = (0.0,1.0,0.045)$.

Before building low-rank representations of $Q$, we reparameterize the problem so that it maps from $[0,1]^{24}$ to the PDE solution value of interest. We do so by using the cumulative distribution function $\Phi$ of a standard Gaussian to define new input variables $\hat{\xi}_i := \Phi(\xi_i)$. Now $(\hat{\xi}_i)_{i=1}^{24}$ are uniformly distributed on the unit hypercube, and the function $Q$ to be approximated is $(\hat{\xi}_i)_{i=1}^{24} \mapsto  u \left ( 0.7, \Phi^{-1}(\hat{\xi}_1),\ldots, \Phi^{-1}(\hat{\xi}_{24}) \right ) $.

\begin{figure}
\begin{center}
\begin{tabular}{m{0.5cm}|D|C|C}
 & $\log(\textrm{squared error})$ &  $\log(\textrm{\# evaluations})$ &  maximum rank \\
\hline
P1 &
\includegraphics[clip,trim=10 10 0 0,scale=0.65]{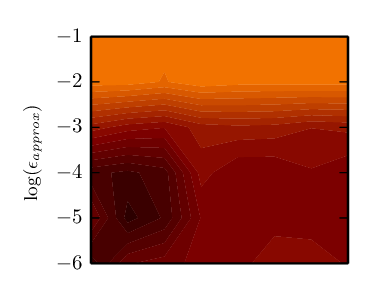} &
\includegraphics[clip,trim=7 10 1 0,scale=0.65]{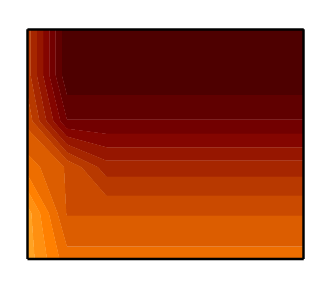} &
\includegraphics[clip,trim=7 10 0 0,scale=0.65]{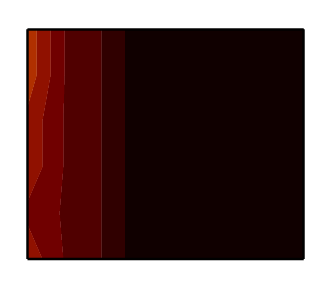} \\
P2 &
\vspace{-2pt}\includegraphics[clip,trim=10 10 0 10,scale=0.65]{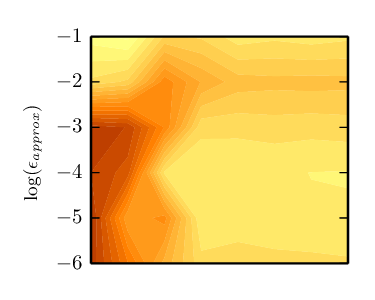} &
\vspace{-2pt}\includegraphics[clip,trim=7 10 0 10,scale=0.65]{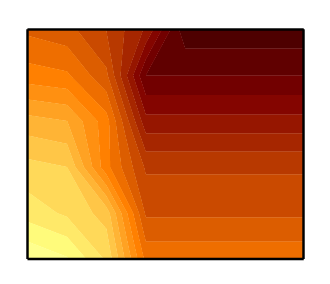} &
\vspace{-2pt}\includegraphics[clip,trim=7 10 0 10,scale=0.65]{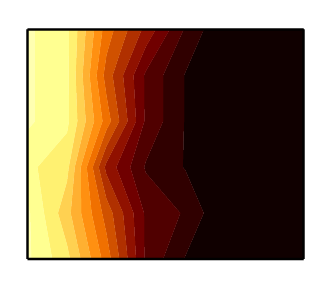} \\
P3 &
\vspace{-4pt}\includegraphics[clip,trim=10 0 0 5,scale=0.31]{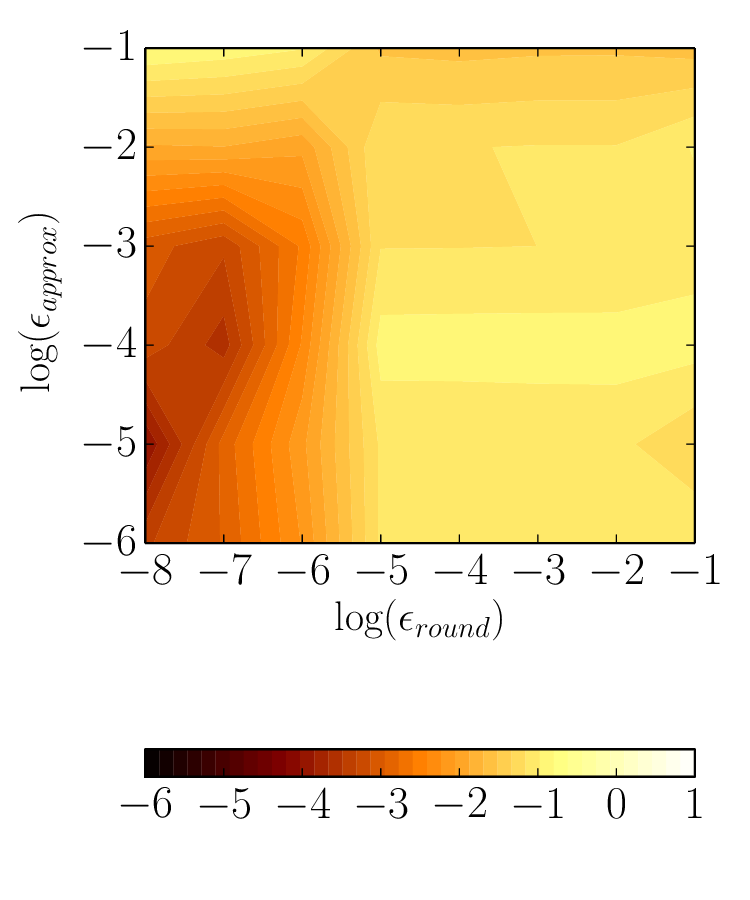} \vspace{-20pt} &
\vspace{1pt}\includegraphics[clip,trim=10 0 0 5,scale=0.30]{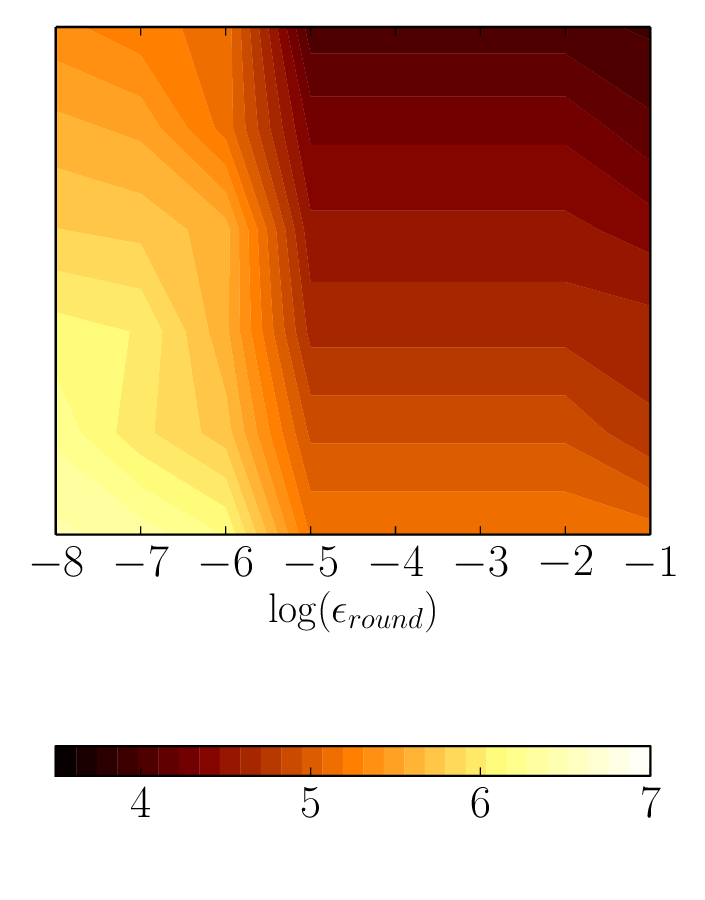} \vspace{-20pt} &
\vspace{1pt}\includegraphics[clip,trim=10 0 0 5,scale=0.30]{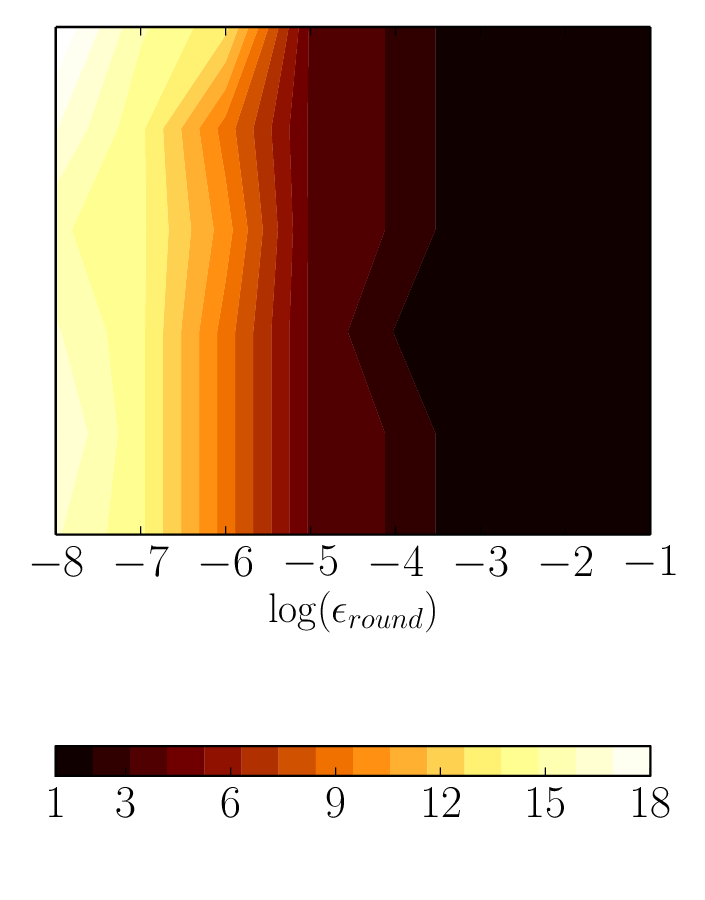}
\vspace{-20pt}
\end{tabular}
\caption{Squared error (left column), number of function evaluations (middle column), and maximum rank (right column) for three different configurations of the elliptic PDE problem, corresponding to different combinations of $(a,\sigma^2,l)$. In particular, the top row corresponds to $(0.0,0.1,0.125)$, middle row corresponds to $(0.5,1.0,0.045)$, and the bottom row corresponds to $(0.0,1.0,0.045)$. Contours in the left two columns represent $\log$-quantities.}
\label{fig:ellipres}
\end{center}
\end{figure}

Our numerical experiments investigate how the final approximation error, number of function evaluations, and maximum rank change with the fiber approximation tolerance $\approxeps$ and the rounding tolerance $\roundeps$. We fix the cross approximation parameters to $\crossdelta=10^{-3}$ and a maximum of 3 sweeps. The rank adaptation parameters are fixed to $\texttt{kickrank}=5$ and a maximum of 4 rank adaptations. We use Legendre polynomials to approximate the fibers and initialize these approximations at degree two. Results for the three problem setups shown in Figure~\ref{fig:ellipres}. 

Several interesting patterns are apparent in the results. First, consider the ranks for each problem, shown in the rightmost column of Figure~\ref{fig:ellipres}. We immediately see that the fiber approximation accuracy has essentially no impact on the maximum rank found by rank adaptation. The maximum rank found by the adaptation procedure only changes as the rounding tolerance decreases. Furthermore, we see that the maximum rank increases as we move down the table. This trend corresponds to the increasing difficulty of each problem, and reflects slower decay of the singular values in decompositions of $u$. Higher ranks for P3 confirm that approximation is more difficult with $a=0$ than with $a=0.5$.  We also note that unlike array-based tensor-train algorithms or the spectral tensor-train algorithm, the maximum rank attainable by the numerical procedure is \textit{not bounded} by discretization level. Rounding is critical to restricting the growth of the rank.

The function evaluation count, shown in the middle column of Figure~\ref{fig:ellipres}, also displays some interesting properties. In particular, for each of the three models, the contours exhibit two separate regimes. Before the rounding tolerance becomes tight enough to cause significant increases in the rank, the number of function evaluations is essentially unaffected by $\roundeps$. The number of evaluations is only affected by the fiber approximation accuracy $\approxeps$. This behavior makes sense because the number of function evaluations is roughly proportional to $\mathcal{O}(ndr^2)$, where $n$ can be thought of as an average number of function evaluations for each fiber---and in this first regime, the rank is constant. Once the rank starts increasing, the number of function evaluations grows with both tighter fiber approximation tolerance and tighter rounding tolerance. The number of evaluations increases more rapidly with tighter rounding tolerance, because reducing the latter produces a rapid increase in rank.

The left column of Figure~\ref{fig:ellipres} shows that the error exhibits a pattern similar to the number of function evaluations. In particular, the error is fairly constant until the rank starts increasing. % \todo{Earlier this said `decreasing.' Please check that I didn't mess up the meaning.}
Once the rank starts increasing, however, then both the fiber approximation tolerance and the rounding tolerance affect the error. Furthermore, if we fix the rounding tolerance and decrease the approximation tolerance, we see a rapid change in approximation error followed by a relatively large plateau area---suggesting that below a given $\approxeps$, either the accuracy of fiber approximations can no longer appreciably increase or further gains in accuracy have little effect on the overall approximation error. Overall, these results suggest that it is possible to find a reasonable value for the rank before increasing the fiber approximation accuracy. Future work will require investigating how to jointly adapt these two parameters in the best way. 

Finally, comparing the error plots for P2 and P3, we observe that the contours look quite similar, but that the ranks and numbers of function evaluations are larger in the bottom row (P3) than in the middle row (P2). This characteristic is highly desirable for an ``adapt-to-tolerance'' scheme, as a given setting for the tolerances yields similar errors but appropriately larger computational effort.

\section{Conclusions}

We have developed new algorithms and data structures for representing and computing with the functional tensor-train decomposition. The algorithms are both locally and globally adaptive, and enable robust approximation, integration, and differentiation schemes. The primary way we achieve adaptation is by following the Chebfun paradigm of continuous computation: we consider continuous extensions of the CUR and QR matrix decompositions that naturally embed adaptivity within the resulting algorithms.

We have demonstrated that such schemes can enable several orders of magnitude more accurate approximation \revise{than the state-of-the-art STT or QTT approaches} for the same computational expense. These savings arise because using continuous linear algebra provides a flexible way of incorporating and exploiting more than just low-rank structure. In particular, our algorithms can exploit the regularity of the target function and adapt to its features. 
Pivots for skeleton approximation are found via continuous optimization. Ranks can be adapted via a rounding procedure. And the FT decomposition offers an enormously flexible approach to univariate fiber approximation: different bases or approximation schemes may be chosen for different input coordinates, or even for parallel fibers in the same dimension; and these approximations can be adapted or refined on a fiber-by-fiber basis. Function evaluations can thus be tailored to local features, without following any tensor product structure. Indeed, the FT scheme employs no \textit{a priori} discretization of the target function and does not require specifying a tensor product set of candidate evaluation locations. This characteristic is particularly important for problems that are sensitive to the choice of discretization.  
The data structure that we have developed stores each univariate function independently, which allows for within-dimension differences in parameterization and increased compression rates when compared with the discrete TT cores.

Future algorithmic extensions will incorporate the \texttt{dmrg-cross} technique of \cite{Savostyanov2011} into a continuous rank-revealing algorithm, which should be more efficient than the rank-adaptation algorithm described here. Future work will also seek to exploit low-rank structure in inference and control. For example, we would like to extend the optimal stochastic control framework described in~\cite{Gorodetsky2018} to problems that are control-affine and avoid discretization of the state space altogether.

\section*{Acknowledgments}
This work was supported by the National Science Foundation through grant IIS-1452019, and by the US Department of Energy, Office of Advanced Scientific Computing Research under award number DE-SC0007099.

\section*{References}
\bibliography{ftrain.bib}

\appendix

\end{document}